\documentclass[<option>]{elsarticle}
\usepackage{amsmath}
\usepackage{pdflscape}
\usepackage{multicol}
\usepackage{color,soul}
\usepackage{amsmath,amssymb}
\usepackage{mathtools}
\usepackage{multirow}
\usepackage{rotating}
\usepackage[table]{xcolor}
\usepackage[a4paper]{geometry}
\usepackage{paralist}
\usepackage{bm}
\usepackage{mathabx}
\usepackage{caption}
\usepackage{float}
\usepackage{graphicx}
\usepackage{caption}
\usepackage{subcaption}
\usepackage{graphicx}
\usepackage{graphics} 
\usepackage{epsfig} 
\usepackage{graphicx} \usepackage{epstopdf}
\usepackage[colorlinks=true]{hyperref}
\hypersetup{urlcolor=blue, citecolor=red}
\begin{document}

\begin{frontmatter}

\title{Oncolytic virotherapy for tumours following a Gompertz growth law}
\author[mymainaddress]{Adrianne L. Jenner}
\author[mymainaddress]{Peter S. Kim}
\author[mysecondaryaddress]{Federico Frascoli}

\address[mymainaddress]{School of Mathematics and Statistics, University of Sydney, Sydney, NSW, Australia}
\address[mysecondaryaddress]{Department of Mathematics, Faculty of Science, Engineering and Technology, Swinburne University of Technology, Melbourne, VIC, Australia}

\begin{abstract}
We present a mathematical model describing oncolytic virotherapy treatment of a tumour that proliferates according to a Gompertz growth function. We present local stability analysis and bifurcation plots for relevant model parameters to investigate the typical dynamical regimes that the model allows. The model shows a singular equilibrium and a number of nonlinear behaviours that have interesting biological consequences, such as long-period oscillations and bistable states where two different outcomes can occur depending on the initial conditions.

Complete tumour eradication appears to be possible only for parameter combinations where viral characteristics match well with the tumour growth rate. Interestingly, the model shows that therapies with a high initial injection or involving a highly infective virus do not universally result in successful strategies for eradication. Further, the use of additional, boosting injection schedules does not always lead to complete eradication. Our framework, instead, suggests that low viral loads can be in some cases more effective than high loads, and that a less resilient virus can help avoid high amplitude oscillations between tumours and virus.

Finally, the model points to a number of interesting findings regarding the role of oscillations and bistable states between a tumour and an oncolytic virus. Strategies for the elimination of such fluctuations depend strongly on the initial viral load and the combination of parameters describing the features of the tumour and virus.  
\end{abstract}

\begin{keyword}
Oncolytic Virotherapy, Bifurcation Theory, Therapies, Dosage
\end{keyword}
\end{frontmatter}

\section{Introduction}

Oncolytic viruses constitute a class of targeted anticancer agents that have unique mechanisms of action compared with other therapies. The premise is to genetically engineer viral particles to selectively replicate in and lyse tumour cells. Over the past decade, hundreds of patients with cancer have been treated in clinical trials with oncolytic viruses \cite{liu2007clinical}. Unfortunately, due to the heterogeneous nature of cancer, success is elusive, and there is a growing need to quantify the dependency of treatment outcome on cancer characteristics.

A number of mathematical models have been constructed to understand the dynamics of proliferation and diffusion of such viruses in cancerous and healthy tissues. Zurakowski and Wodarz \cite{zurakowski2007model} developed a mathematical model for the \textit{in vitro} behaviour of the oncolytic virus ONYX-015 infecting cancer cells. The virus, administered in conjunction with a drug that up-regulates a tumour cell's ability to intake viral particles, shows two distinct behaviours: a locally stable steady state and a nonlinear periodic cycle between viral particles and tumour cells, in agreement with dynamics that have been experimentally observed. Their model has also allowed them to suggest strategies that alter the amplitude of oscillations and drive tumour size to low levels. Bajzer \textit{et al.} \cite{bajzer2008modeling} introduced a mathematical model that also exhibits stable oscillations between virus and tumour. Crivelli \textit{et al.} \cite{crivelli2012mathematical} instead derived a cycle-specific age-structured model for virotherapy where cell-cycle-specific activity of viruses has been investigated. Through analysis and simulation of the model, the authors have described how varying minimum cycling time and aspects of viral dynamics may lead to complete eradication of the tumour. Other authors have also focused their modelling on the delay occurring between the initial virus infection of tumour cells and the second wave of infections when viruses burst \cite{cassidy2018mathematical,jenner2018heter}. This phenomenon has also been accounted for by a delay differential equation model for the lapse in the second generation viral infection \cite{kim2018hopf}. 

From the experimental point view, particularly relevant to the present work are the findings by Kim \textit{et al.} \cite{KimPH2011}. These authors have developed an oncolytic virus, modified with immunogenic polymer polyethylene-glycol (PEG) and monoclonal antibody Herceptin, which has exhibited potent anti-tumour behaviour in murine models. While the treatment appears unable to eliminate the tumour completely, it significantly reduces the growth rate of cancer cells with respect to the case of an untreated tumour. Aiming to provide an explanation as to why this treatment seems to fully eradicate a tumour, we have recently proposed a mathematical model calibrated to these experimental results \cite{jenner2018mathematical}. Analysis of common treatment protocols has unveiled some drawbacks of existing strategies and suggested optimal scheduling to maximise therapeutic benefits. 

In this work, we conduct a further investigation of this mathematical model with a particular emphasis on the biologically relevant parameters that control virus-tumour dynamics. After a discussion of the basis for our approach, we introduce a non-dimensional version of the system, which allows us to conduct local stability analysis and analytically determine which parameters can lead to incomplete tumour eradication, as often observed in experimental settings. Then, we present a bifurcation analysis of the model, leading to some nontrivial and, in some case, counterintuitive findings about the viral characteristics that drive a complete tumour eradication. By examining perturbations of the viral dosage and their effect on different dynamical regions, we show how we can achieve a better application of treatment protocols. In particular, the dynamical states displayed by the model when therapies are administered strongly affect the final outcomes. Finally, a discussion of the advantages and limitations of our approach concludes the paper.

\section{Model development}

The dynamics of a tumour treatment administered via a PEG-modified adenovirus conjugated with Herceptin can be captured using a minimal mathematical frameworks as explained in Ref.~\cite{jenner2018mathematical}. Assuming that the immune response is negligible and does not need to be incorporated in the equations, the following model can be proposed:
\begin{align} 
\frac{dU}{d\tau} & = r\ln \left(\frac{K}{U}\right)U -\frac{\beta U\hat{V}}{U+I},\label{AEqs1}\\
\frac{dI}{d\tau} & = \frac{\beta U\hat{V}}{U+I}-d_II,\\
\frac{d\hat{V}}{d\tau} & = -d_V\hat{V}+\alpha d_I I, \label{AEqs4}
\end{align} 
\noindent where $\tau$ is time, $\hat{V}$ represents the density of virus particles at the tumour site, $U$ is the density of susceptible but virus-uninfected tumour cells, $I$ is the density of virus-infected tumour cells and the term $U+I$ corresponds to the total tumour cell population. 

Tumour growth is controlled by nutrients and spatial limitations and is described by a Gompertz function, i.e. $g(U) = r\ln(K/U)U$. Here, $K$ represents the carrying capacity of the tumour and $r$ is its growth rate. This type of expression is well-known to reproduce the experimentally observed evolution of a number of proliferating tumours quite accurately~\cite{laird1964dynamics}. In our framework, the likelihood of a virus infecting a tumour cell is assumed to depend on the number of tumour cells available to infect. To model this, a frequency-dependent function, rather than a simple mass-action term, is introduced: virus particles at the tumour site infect susceptible tumour cells according to the expression $\beta U\hat{V}/(U+I)$, where $\beta$ is the infectivity rate. This also differentiates this model from existing, well-known virus dynamics models, such as systems used for influenza or HIV modelling~\cite{DeLeenheer20031313, Wang200644}.

\begin{figure}[h!]
    \centering
    \includegraphics[width=80mm]{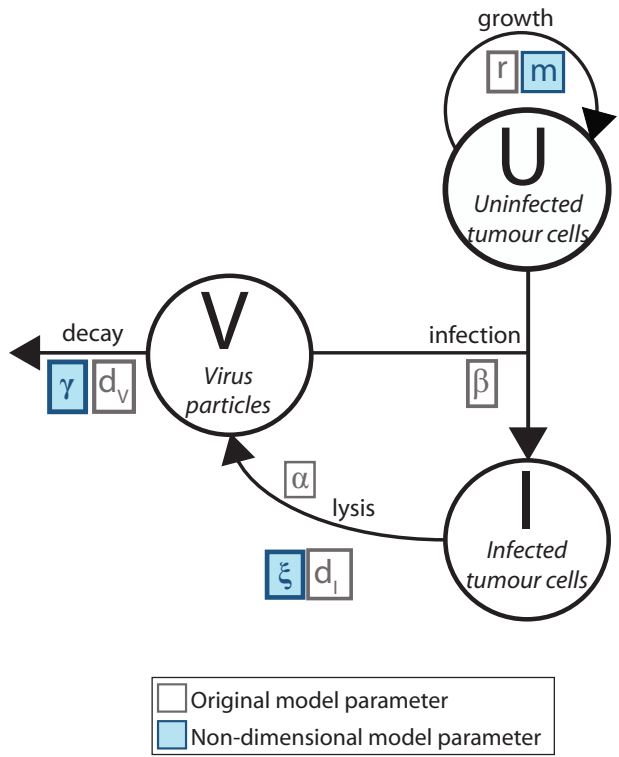}
    \caption{Flow diagram for the interaction between a population of uninfected tumour cells, $U$; virus-infected tumour cells, $I$; and virus particles, $V$. The diagram lists parameters relating to the original model Eqs.~(\ref{AEqs1})-(\ref{AEqs4}), in grey boxes and parameters relating to the non-dimensional form of the model, Eqs.~(\ref{E4})-(\ref{E6}), in blue boxes.}
    \label{Schem11}
\end{figure}
After initial injection, any virus subsequently produced via replication within tumour cells will not have the PEG or Herceptin modifications. To account for this, only single average infectivity and decay rates $\beta$ and $d_V$ are used for the combined populations of original and replicated virus, noting that the population is dominated by naked (replicated) virus over the majority of the time course of the experiments. Fig.~\ref{Schem11} depicts the flow diagram of the three populations described in the equations, and we refer to the original study \cite{jenner2018mathematical} for a discussion on biologically relevant ranges of values of the parameters. 

To proceed with our mathematical analysis, an appropriate change of variables detailed in Appendix A is used to scale the above equations into dimensionless form. The final result is given as follows:
\begin{align}
\frac{dU}{dt} &= m\ln\left(\frac{K}{U}\right)U-\frac{UV}{U+I},\label{E4}\\
\frac{dI}{dt} &= \frac{UV}{U+I}-\xi I,\\
\frac{dV}{dt} &= -\gamma V+\xi I\label{E6}
\end{align}
\noindent where $m = \displaystyle\frac{r}{\beta}, \xi = \displaystyle\frac{d_I}{\beta}, \gamma =\displaystyle \frac{d_V}{\beta}$ and $\hat{\beta} = \displaystyle\beta \alpha$ are dimensionless parameters, and $t$ represents a dimensionless ``time''. This model still follows the schematic given in Fig.~\ref{Schem11} and is the object of the present study. The three parameters $m, \xi$ and $\gamma$ that regulate the behaviour of the system represent tumour growth rate, viral death rate and viral potency (or infectivity), respectively. As a result of the non-dimensionalisation process, where parameters are all scaled by the infectivity rate (see \ref{sec:appendix}), the rate of conversion of uninfected cells $U$ to infected cells $I$ due to the viral load $V$, i.e. the term $\pm\displaystyle \frac{UV}{U+I}$, is not affected by any parameter.

\section{Local stability analysis}

A local stability analysis of Eqs.~(\ref{E4})-(\ref{E6}) shows a number of interesting features. Of particular relevance is the existence of a stable equilibrium corresponding to eradication, which is characterised by a singular Jacobian matrix. This solution can coexist with other equilibria, for example a stable spiral or a stable node, which instead corresponds to incomplete eradication of the tumour. As we will show shortly, this occurrence can give rise to bistability for some biologically relevant parameter ranges. 

\subsection{Equilibrium solutions}

Setting the right-hand-side of Eqs.~(\ref{E4})-(\ref{E6}) to zero, three equilibria are found: (a) a solution at a value for the uninfected cells equalling the carrying capacity, indicating a treatment with no effect; (b) a non-zero solution representing incomplete eradication, characterised by a quiescent tumour despite the viral load being constant and non zero; and (c) an equilibrium at the origin corresponding to complete eradication of the tumour. The populations corresponding to such cases are

\begin{align*}
(a)& \quad U= K, \ \  I=0, \ \ \ V=0; \\
(b)& \quad U= K\exp\left(\displaystyle\frac{\xi}{m\gamma}(\gamma-1) \right)=U^*, \ \ I= \frac{K}{\gamma}(1-\gamma)\exp\left(\displaystyle\frac{\xi}{m\gamma}(\gamma-1) \right) =I^*, \\ 
&\quad V = \frac{K\xi}{\gamma^2}(1-\gamma)\exp\left(\displaystyle\frac{\xi}{m\gamma}(\gamma-1) \right)= V^*;\\
(c)& \quad U = 0, \ \ \ I=0, \ \ \ V=0.\\~\label{nonzeroequil}
\end{align*} 

\noindent The Jacobian of the system is given by

\begin{align}
J = \left(\begin{array}{ccc}
m\ln\left(\displaystyle\frac{K}{U}\right)-m-\displaystyle\frac{VI}{(U+I)^2} & \displaystyle\frac{UV}{(U+I)^2} & -\displaystyle\frac{U}{(U+I)}\\[8pt]
\displaystyle\frac{VI}{(U+I)^2} & -\xi-\displaystyle\frac{UV}{(U+I)^2} &\displaystyle\frac{U}{U+I}\\[8pt]
0 & \xi& -\gamma
\end{array}\right),
\end{align}

and we discuss the character of the eigenvalues for the above equilibria below.

\subsection{Stability of ineffective treatment equilibrium: $U=K$, $I=0$, $V=0$}

The first equilibrium (a) corresponds to a failed treatment where uninfected tumour cells $U$ grow to the system's carrying capacity $K$ and no viral particle survives. Evaluating the Jacobian at this point gives
\begin{align*}
J = \left(\begin{array}{ccc}
-m & 0 & -1\\
0 & -\xi &1\\
0 & \xi& -\gamma
\end{array}\right),
\end{align*}
which gives rise to the following characteristic equation:
\begin{equation}\label{ch-eq-partial}
\rho(\lambda; m, \gamma, \xi) = -(\lambda+m)\left(\lambda^2+(\xi+\gamma)\lambda+\xi(\gamma-1)\right).
\end{equation}
The overall stability of this equilibrium depends on the roots $\lambda_2$ and $\lambda_3$ of the quadratic factor, because the root $\lambda_1 = -m$ of the linear factor is negative, since the growth rate $m>0$. After calculating $\lambda_2$ and $\lambda_3$, we find that the equilibrium is either a stable node or stable focus when $\xi+\gamma>0$ and $\xi(\gamma-1)>0$. Since the parameter values in this model are considered to be always positive, the first condition holds. The second condition implies that, if $\gamma<1$, the equilibrium is unstable, and vice versa for $\gamma>1$. As we will show shortly, a one-parameter continuation in $\gamma$ shows that at $\gamma = 1$ a branch point is present and the treatment is always ineffective for a decay rate $\gamma>1$. Intuitively, if the virus dies too quickly, no infection can occur.

\begin{figure}[h!]
\centering
\includegraphics[width=0.44\textwidth]{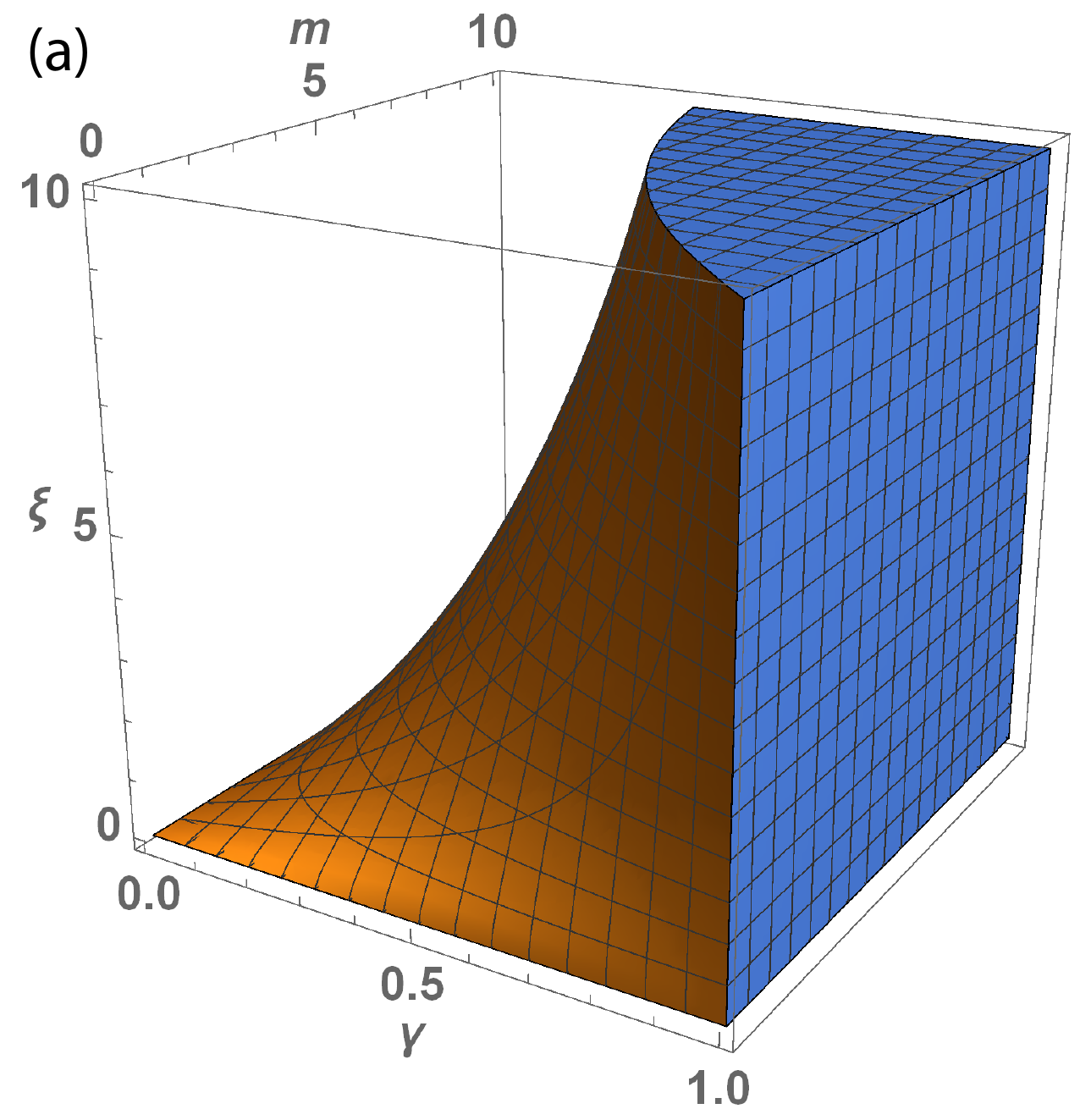}
\includegraphics[width=0.44\textwidth]{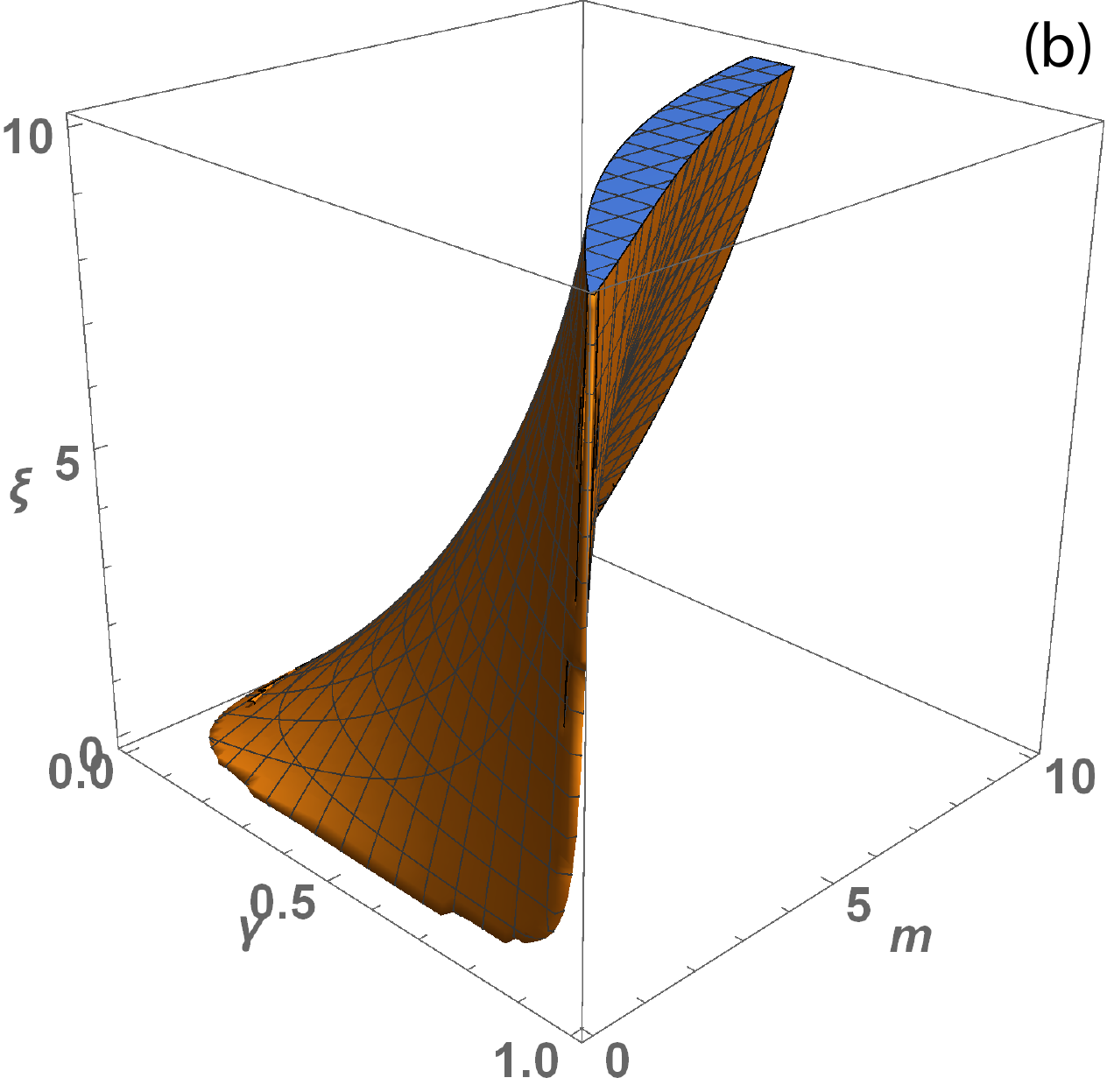}\\[5pt]
\includegraphics[width=0.44\textwidth]{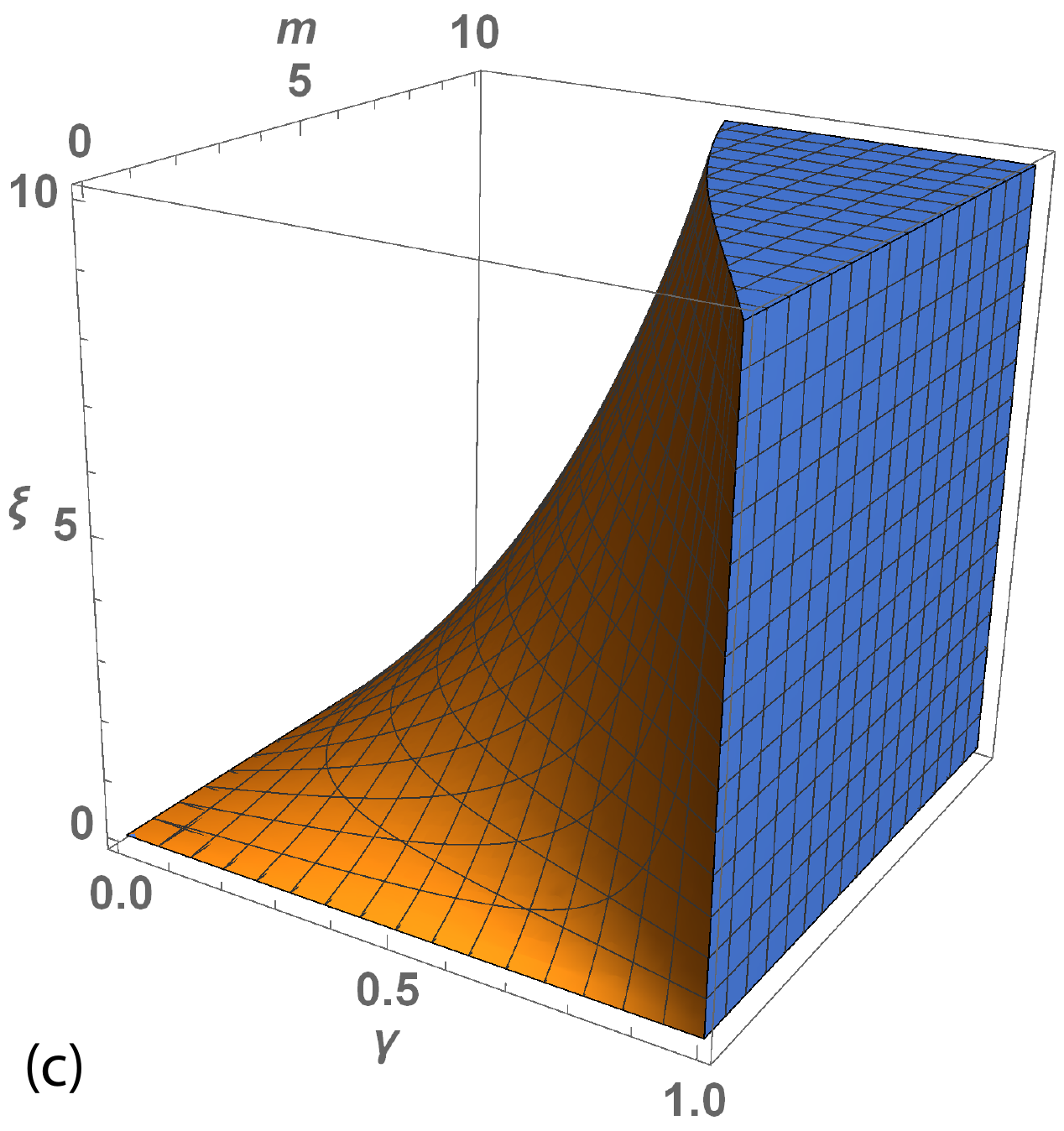}
\includegraphics[width=0.44\textwidth]{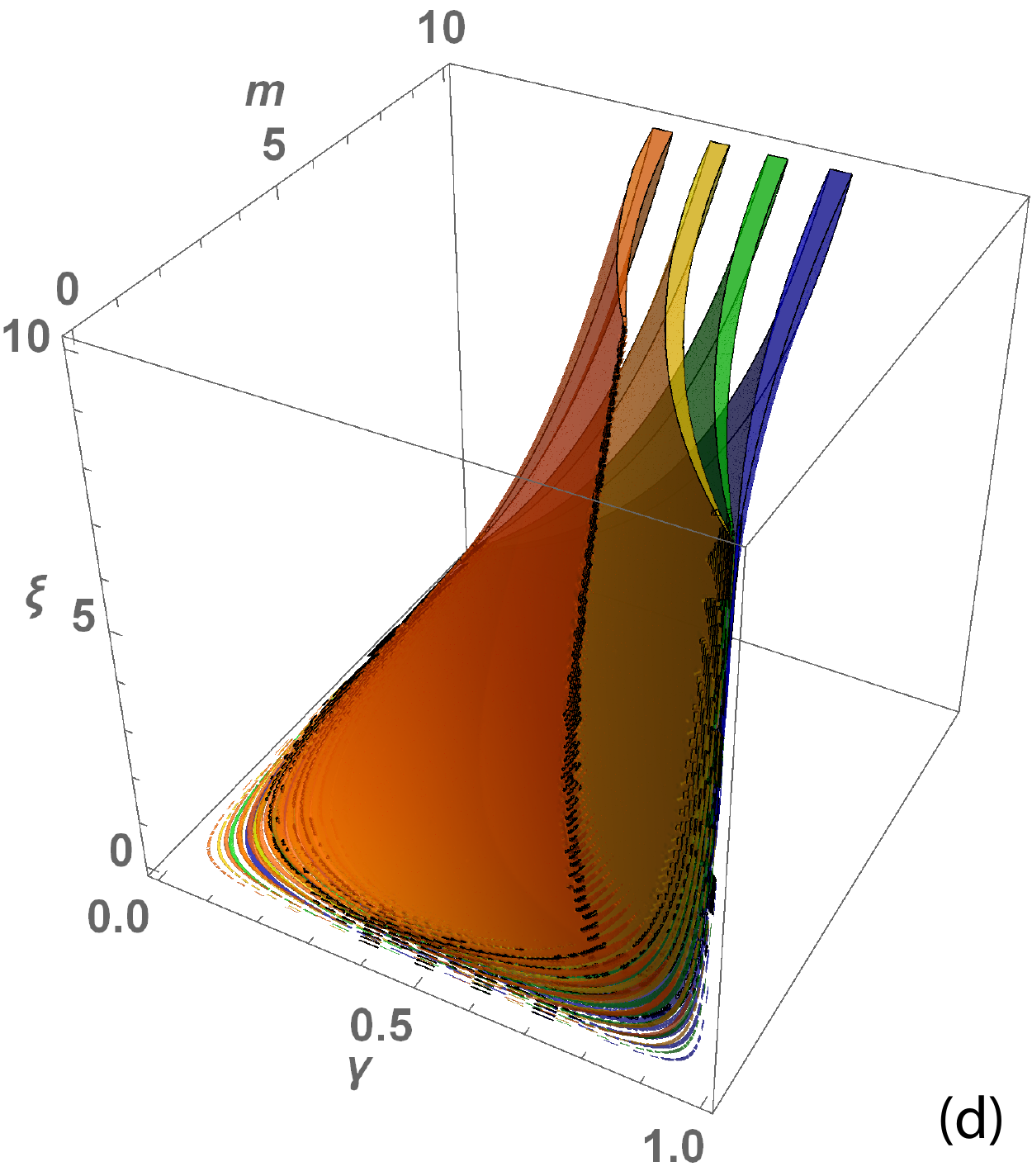}
\caption{Regions representing the stability of the nonzero equilibrium, (a)-(c), and the influence of system parameters on tumour cell numbers at the equilibrium value $U^*$. In (a), the section of parameter space where the non-zero equilibrium is stable is shown. Note that (b) represents the volume in $(\xi,m,\gamma)$ giving rise to a stable node solution for the equilibrium $(U^*, I^*, V^*)$, whereas (c) is the section for a stable spiral. Combining the regions in (b) and (c) gives the volume in (a). Plot (d) is the stable parameter space for different values of $U^*$, within the following intervals: orange for $0.2<U^*<0.25$, yellow for $0.35<U^*<0.4$, green for $0.5<U^* < 0.55$ and blue for $0.65<U^*<0.7$. Note that these ``slices'' are almost symmetrical.}\label{nsregion}
\end{figure}

\subsection{Stability of partial eradication solution: $U=U^*$, $I=I^*$, $V=V^*$}

The model admits a second, non-zero equilibrium where a small tumour mass coexists with virus particles. The characteristic equation for this solution, after substituting $U^*,I^*,V^*$ in the Jacobian, is given by 
\begin{equation}\label{characeqn}
\rho(\lambda) = -\lambda^3-\lambda^2(\gamma+m+\xi)+\lambda\left(\gamma m(\xi-1)+\frac{\xi^2}{\gamma}-\xi(2m+\xi)\right)+\gamma m \xi (\gamma-1).
\end{equation}
\noindent For this cubic, the Routh-Hurwitz criterion is used to deduce the parameter values that produce three roots with negative real part. This criterion states that, given a general cubic of the form  $\rho(\lambda)=a_0\lambda^3+a_1\lambda^2+a_2\lambda+a_3$, two conditions need to be met simultaneously for all roots to have negative real parts, i.e.

 \begin{equation*} (a)~\frac{a_1a_2-a_0a_3}{a_1}<0 \qquad \text{and} \qquad (b)~a_3<0
\end{equation*}
with, in our case, $a_0 = -1$, $a_1=-(\gamma+m+\xi)$, $a_2= \left(\gamma m(\xi-1)+\frac{\xi^2}{\gamma}-\xi(2m+\xi)\right)$ and $a_3=\gamma m \xi (\gamma-1)$.
Condition (b) is easily satisfied for $0<\gamma<1$, given that all parameters are assumed to be positive. Condition (a) requires that $a_1a_2>a_0a_3$, since $a_1<0$. The region in the $\xi,m,\gamma$ parameter space that satisfies this condition can be numerically computed and is depicted in Fig.~\ref{nsregion}(a). Using the discriminant of Eq.~(\ref{characeqn}) and imposing the appropriate conditions, subsections of that region corresponding to a stable node or stable spiral are illustrated in Fig~\ref{nsregion}(b, c). Note that all regions are smooth and connected. 

It is also interesting to consider which parameter regimes result in a low tumour burden (or threshold) $U_T$. To visualise how the value of the equilibrium $U^*$ changes as a function of parameter values, we can compute the regions of parameters space satisfying the following equality for a given threshold $U_T$:

\begin{equation}
\xi =  \frac{m}{\gamma-1}\ln\left(\frac{U_T}{K}\right).\label{Ustarcontoursurfaces}
\end{equation}

\noindent Plots for four different $U_T$, varying within intervals, are shown in Fig.~\ref{nsregion}(d).  The regions are roughly symmetric, with parameter $\gamma$ being the major contributor to changes in $U^*$ values. For example, when $\gamma \lessapprox 0.5$, there is a set of $\xi$ and $m$ values resulting in $0.20\lessapprox U^*\lessapprox 0.25$. Since $m$ represents the growth rate of tumours and $U^*$ is almost insensitive to its variations, our analysis indicates that a value of $\xi$ can always be chosen to decrease the volume of the tumour, as long as the decay rate $\gamma$ is low (i.e. the virus does not decay too quickly).

\subsubsection{Stability of full eradication solution: $U=0$, $I=0$, $V=0$}

The last equilibrium of the model represents the case of complete eradication, where all variables are zero. As anticipated, the Jacobian is singular due to the presence of logarithmic and rational terms in $U$ and $(U+I)$ respectively. An analytical treatment is not possible and, in particular, the presence of logarithmic terms $m\ln(K/U)$ or its source in Eq.~(\ref{E4}), i.e. $m\ln(K/U)U$, is not treatable with straightforward expansions for $U\to 0$. A different approach based on numerical integration and computation of eigenvalues under specific assumptions on $U$, $I$ and $V$ is instead used and will be discussed in detail in the next section. 

As far as the equilibrium's stability is concerned, it turns out that the eradication solution can be stable or unstable, depending on the value of model parameters. As a general rule, we observe that parameter sets where $\xi$ is high, corresponding to a potent viral load, tend to yield a stable equilibrium as long as the decay rate $\gamma$ is not excessive. This suggests that the engineered virus has to be potent and sufficiently resilient: one characteristic alone is not sufficient. If, for example, the virus has potency $\xi$ but dies too fast, then the equilibrium turns into an unstable point and no eradication is possible. A clear picture of how eradication depends on viral characteristics will emerge with the aid of bifurcation plots, which are discussed in the next section. 

\section{Characteristic dynamical regimes} 
\label{section:4}
\begin{figure}[h!]
 \centering
 \includegraphics[width=0.48\textwidth]{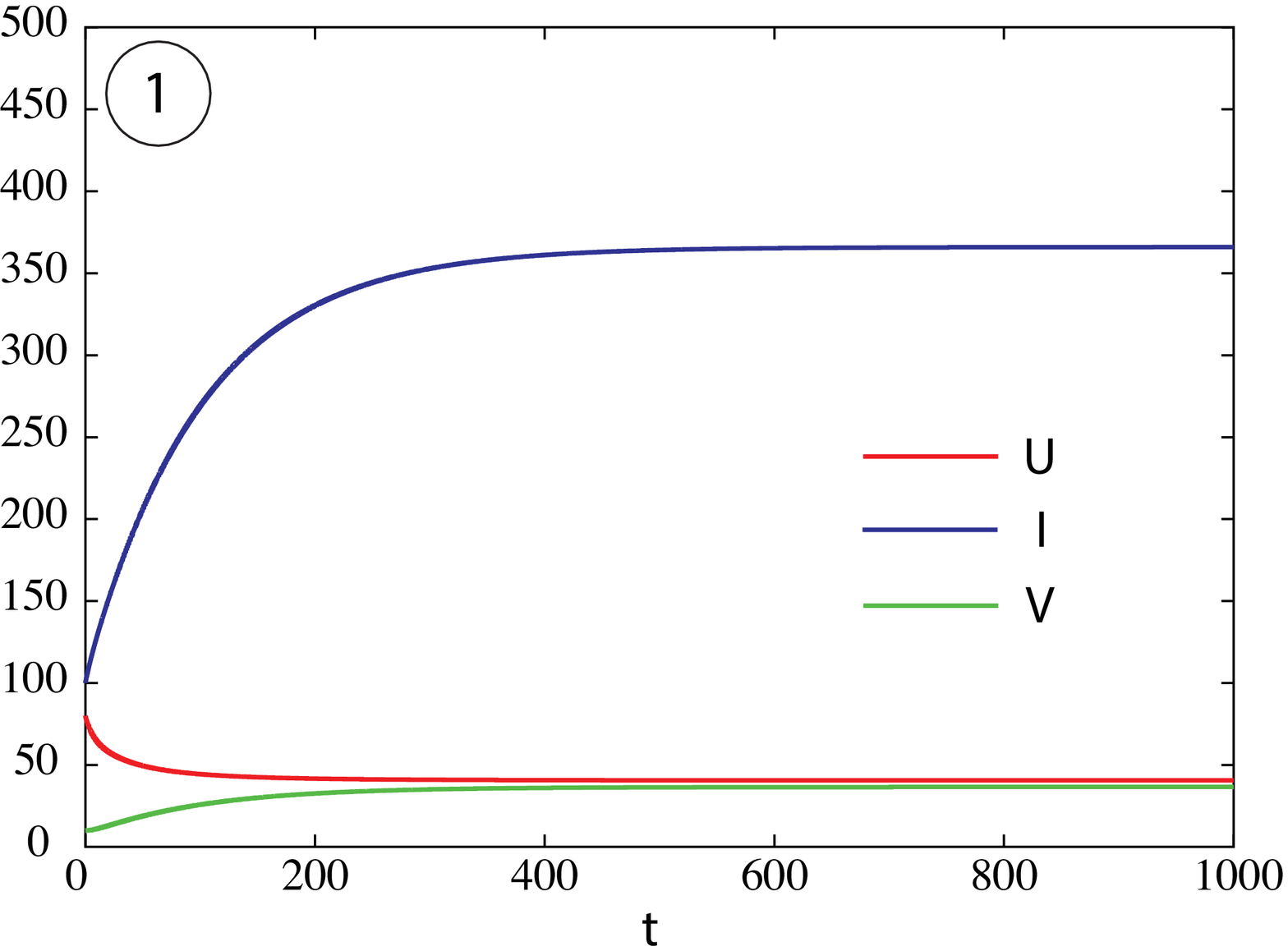}
 \includegraphics[width=0.48\textwidth]{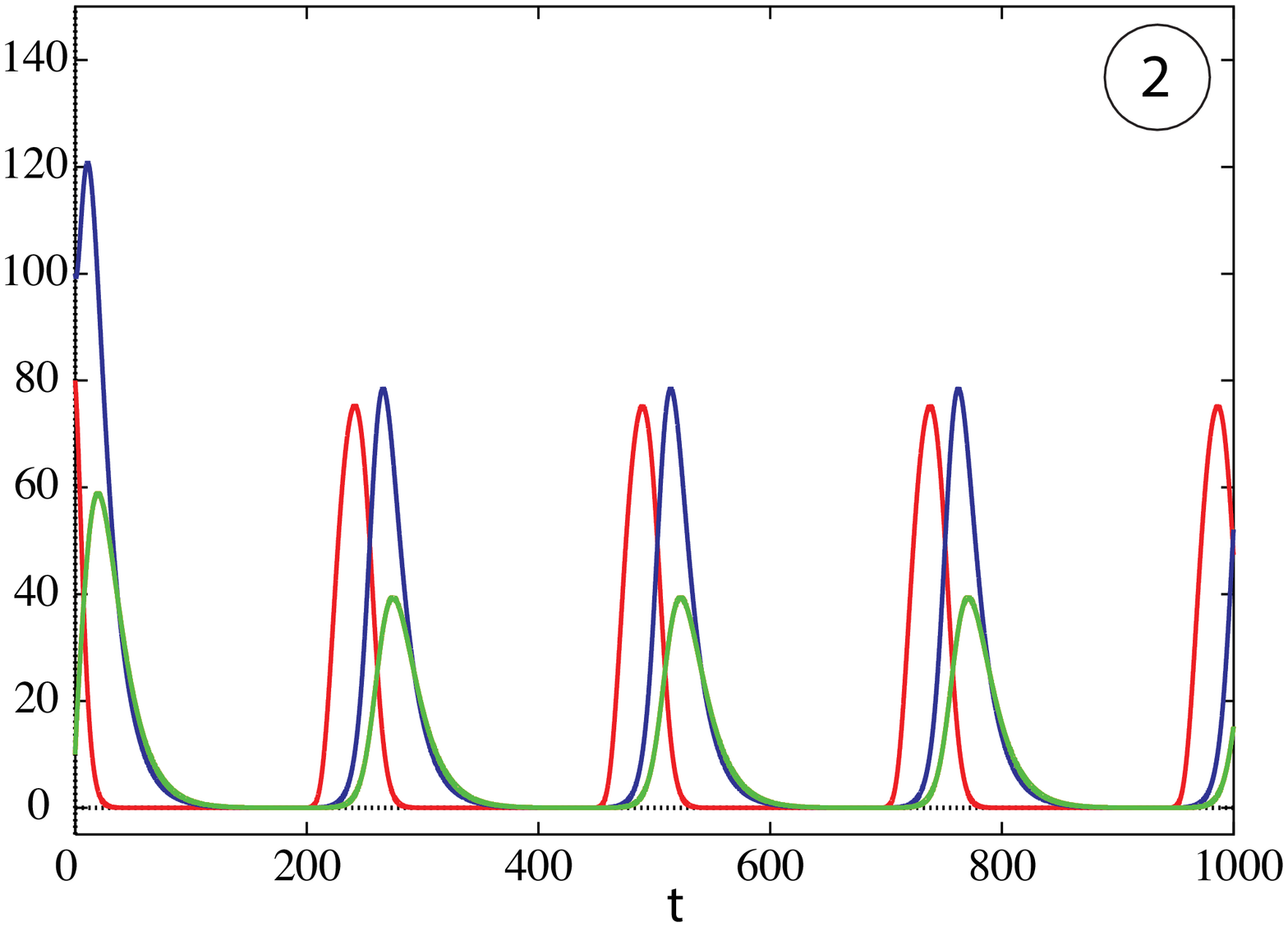}\\[6pt]
 \includegraphics[width=0.48\textwidth]{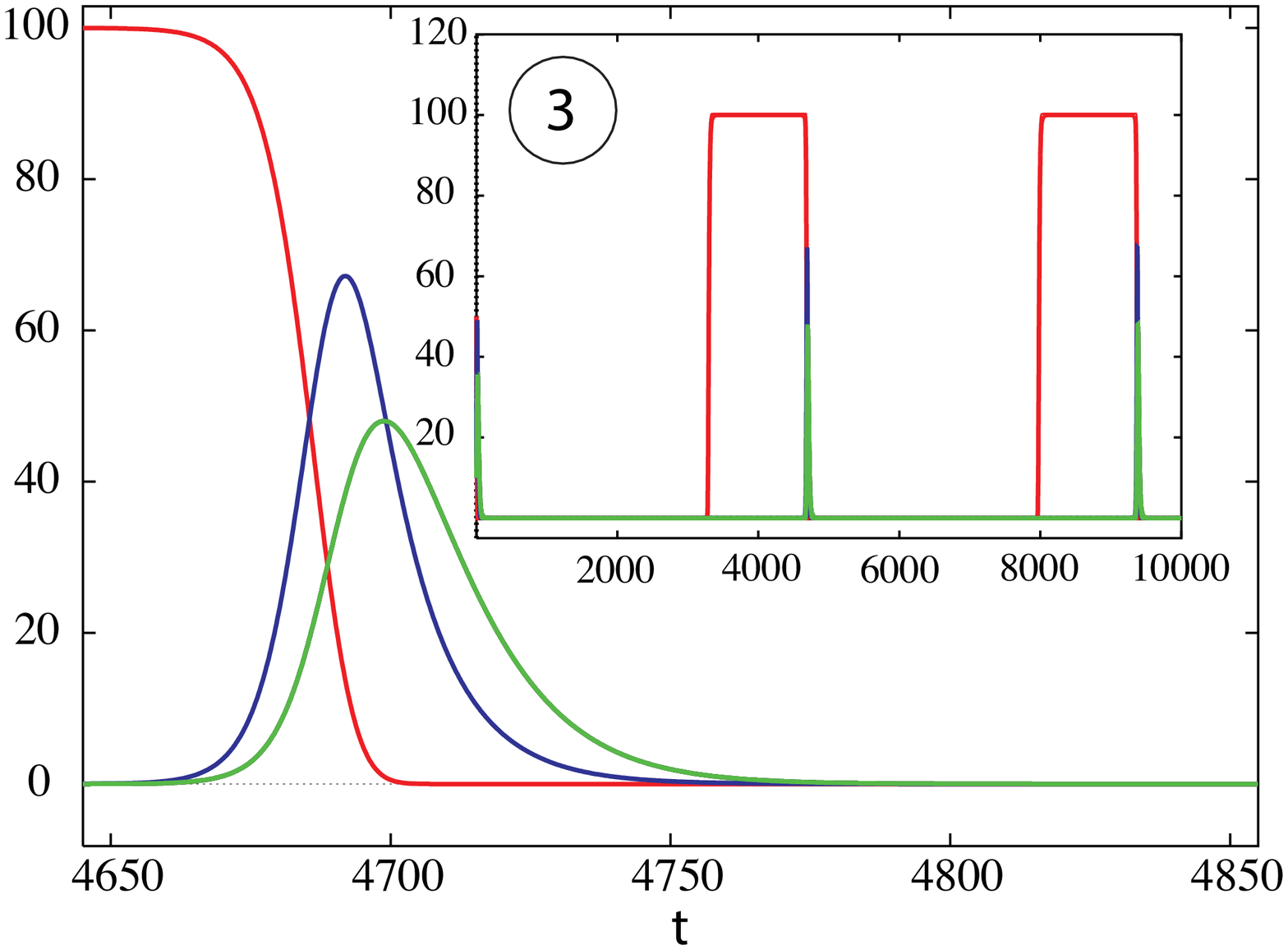}
 \includegraphics[width=0.48\textwidth]{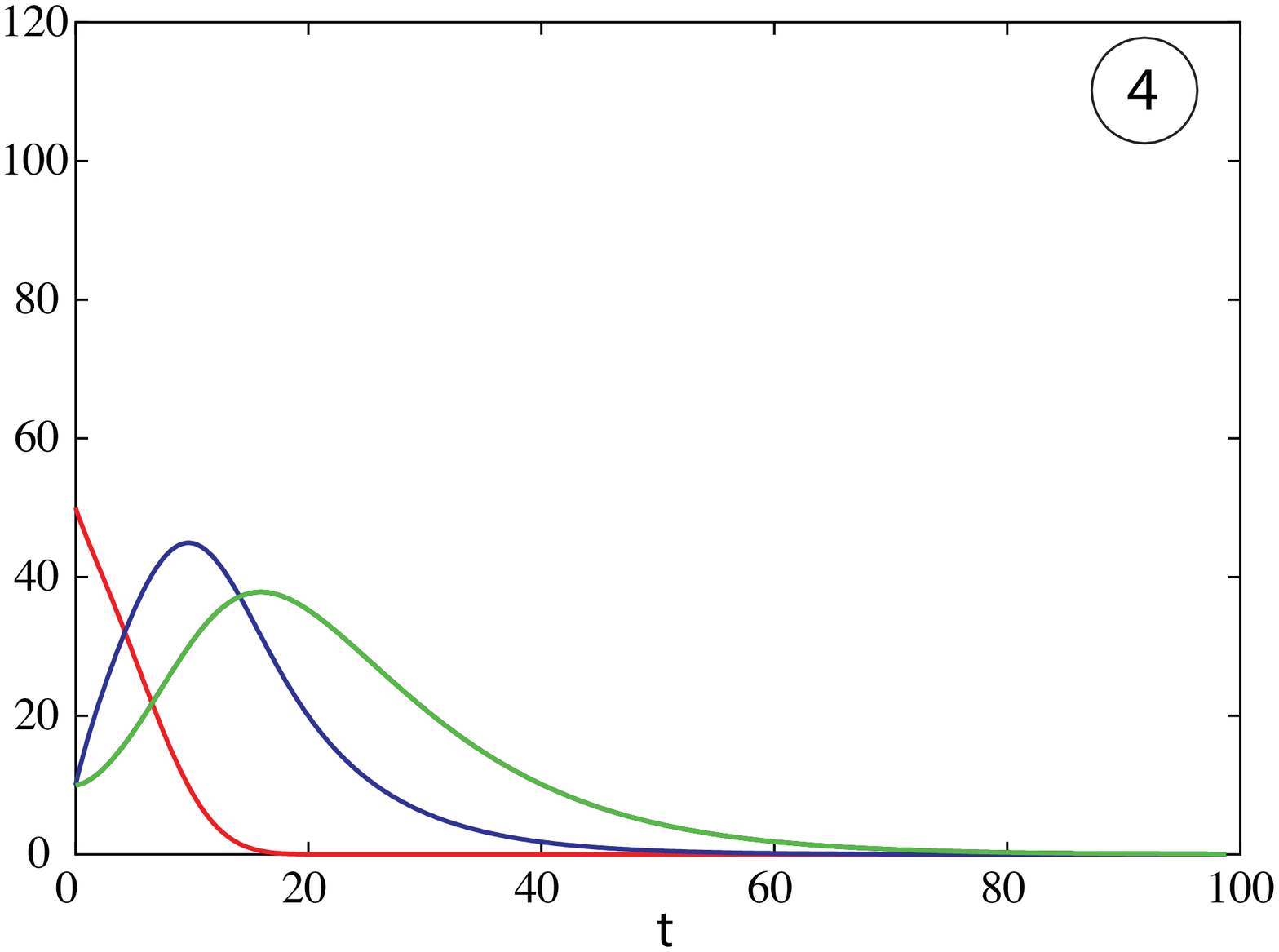}
 \caption{Numerical simulations of Eqs.~(\ref{E4})-(\ref{E6}) demonstrating different types of dynamics, for initial conditions $U(0)=50, I(0) = 10, V(0)= 10$ and fixed parameters $m = 0.1, \gamma = 0.1$, whereas values for $\xi$ are increasing from case (1) to (4). Type (1) corresponds to a stable co-existence of virus and tumour due to incomplete eradication, occurring at $\xi = 0.01$, (2)  depicts a stable oscillatory solution for $\xi = 0.06$, (3) shows stable long period oscillations of almost ``square wave'' shape for $\xi = 0.097$ and (4) is a case of complete eradication for $\xi=0.12$. Note that the carrying capacity is chosen as $K=100$.}\label{Fdyn}
\end{figure}

The model supports a number of dynamical regimes that are interesting both from the biological and mathematical point of view. In Fig.~\ref{Fdyn}, four distinctive behaviours associated with the equilibria previously described are presented. Case (1) is an example of an equilibrium solution where the virus co-exists with uninfected and infected tumour cells, i.e. the case $U=U^*$, $V=V^*$ and $I=I^*$. The time series is for an attracting node, but similar long-term dynamics exists for the case of an attracting spiral, with the only difference being an initial, oscillatory transient that then damps down to a plateau. Note how the uninfected cells $U$ are the first to reach the equilibrium $U^* = K\exp(\frac{\xi}{m\gamma}(\gamma-1))$, which corresponds, for the chosen parameters, to $U^* \approx 40.65$. Case (2) corresponds to stable oscillations, characterised also by an almost quiescent phase where the system variables are close to zero and periods of growth and decay of cells and virus. Generally, we observe that this ``refractory'' state tends to have a longer duration than the active phase. Also in this case, the uninfected cells $U$ are the first to grow, with a subsequent increase in the infected cells $I$ and then in the virus load $V$. As we will see shortly with a bifurcation analysis, the duration of the ``rest'' and ``active'' phases of oscillations depends on the system parameters and changes continuously from case (2) to the limiting case (3). This is an extreme scenario where the system oscillates between two long plateaus of quasi-complete eradication (i.e. $U=I=V\approx 0$) and  quasi-ineffective treatment (i.e. $U\approx K=100, ~I=V\approx 0$). The inset shows the almost square-wave appearance of the system's trajectories on a long time scale, whereas the switch from the two states is illustrated in the main figure, showing how the growth in $I$ and $V$ causes the uninfected cell numbers to decrease. It is important to note that the system cannot stabilise on either equilibra, because for the parameters chosen and as it will be evident shortly from bifurcation results, both equilibria are unstable. 

Finally, a complete eradication solution is depicted in case (4). Although for the chosen initial conditions and parameters the model shows a monotonic decline to zero for $U$, other examples have been found where $U$ first shows a maximum, followed by an exponential decrease. Also in this final case, as for the other three scenarios just discussed, we observe that $U$ is the fastest to reach its equilibrium value, with $I$ and $V$ following. 

\begin{figure}[h!]
\centering
\includegraphics[width=0.49\textwidth]{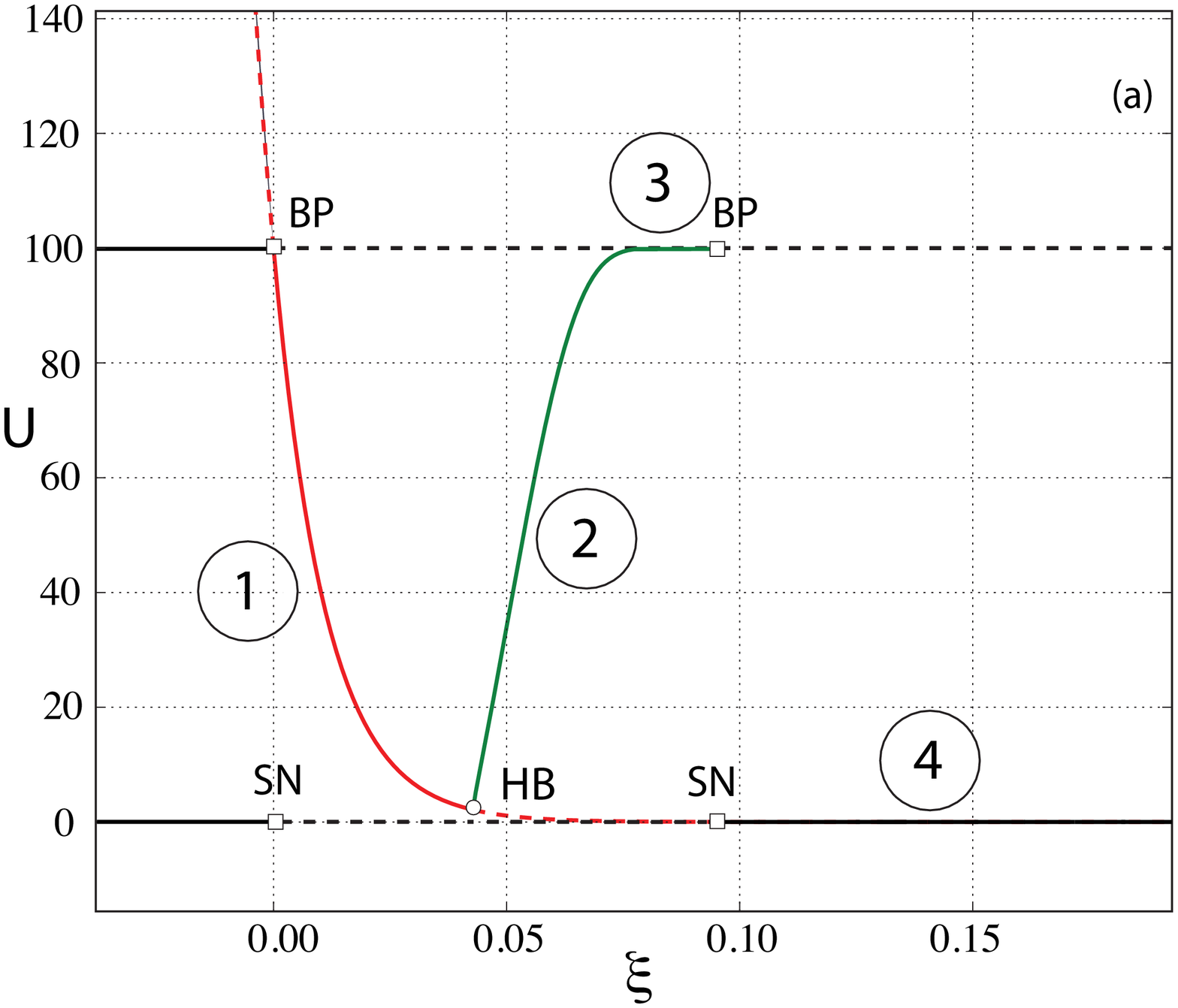}
\includegraphics[width=0.49\textwidth]{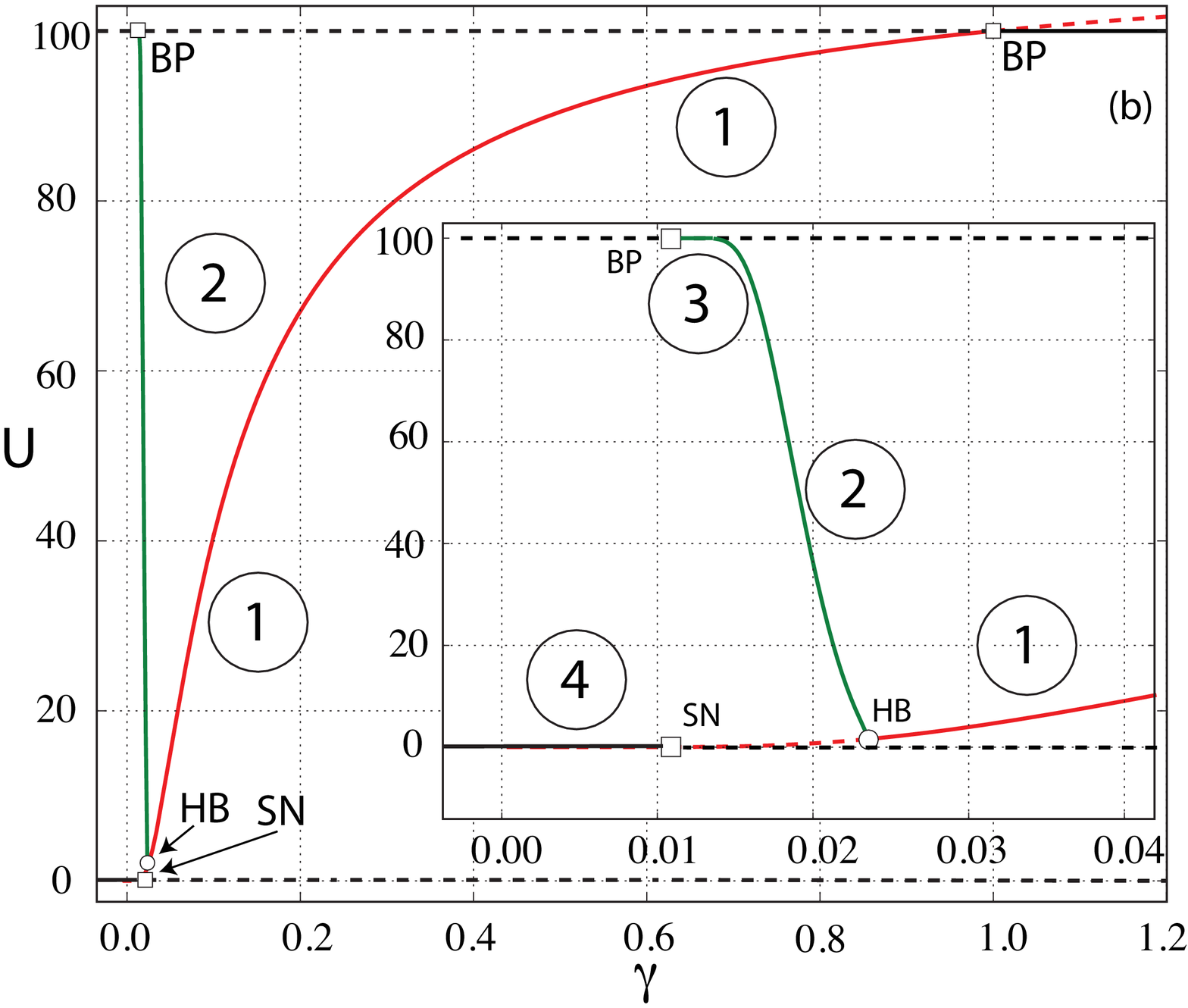}
\caption{Examples of typical bifurcation plots in one parameter for the model, with the case of $\xi$ in (a) and $\gamma$ in (b), both versus $U$. Numbers correspond to the dynamical regimes illustrated in Fig.~\ref{Fdyn} and, for the case of periodic orbits originating from a Hopf bifurcation, only the maximum value of $U$ is shown. For (a), the other model parameter are $m=0.1$, $\gamma = 0.1$. Note that the switch to regime $4$ (complete eradication) occurs when the branch of periodic orbits (in green) ceases to exist, for a value $\xi \approx 0.098$. Similar results for a continuation in $\gamma$ are shown in (b), with the switch to case $(4)$ dynamics also occurring in correspondance of a branch point for the periodic orbit, at $\gamma \approx 0.0103$. An inset with a magnification on the area that shows the richest dynamical variability is also shown. The value of the other, fixed parameters are given in this case by $m=0.1$, $\xi = 0.01$. In both cases, solutions for negative $\xi$ and $\gamma$ have been included for reasons of consistency, but do not correspond to any biologically meaningful state.}\label{Fbif}
\end{figure}

To appreciate where these regimes occur and how the parameters influence their existence, two bifurcation plots with respect to system variables $\xi$ and $\gamma$ versus $U$ are presented in Fig.~\ref{Fbif}. In both plots, stable branches are indicated with continuous lines, whereas unstable ones are dashed. The two black branches at $U=0$ and $U=K=100$ indicate the full eradication and failed treatment solutions, respectively. The red line indicates the partial eradication case, where a non-zero value for the tumour volume and the viral load is present. Numbers point to areas where the typical dynamics just discussed can be found. 

For the case of a codimension one plot with respect to $\xi$ (Fig.~\ref{Fbif}(a)), two branch points are present: one at $U=100$ and $\xi=0$, where the partial eradication solution coalesces with the failed treatment case, and a second at $U=100$ and $\xi \approx 0.098$ where the oscillatory, stable branch (green line) terminates. This branch originates from a supercritical Hopf bifurcation (HB), which causes the initial partial eradication branch to lose its stability. Note how, at this value of $\xi$, a change in the stability of the eradication solution $U=0$ (black line) also happens, with a saddle-node bifurcation (SN) occurring and a stable, fully eradicating regime appearing for $\xi>\xi_{SN} \approx 0.098$. This eradication solution branch regains its stability at $\xi =0$ through another saddle-node bifurcation (SN). Note also that the partial and full eradication branches (i.e. red and black lines, respectively) do not intersect. Finally, let us remind to the reader that solutions for parameter values that are negative do not bear any biological value.

It is worth noting that, along the red branch of coexisting solutions, $U$ can span a large range of values, with $U$ increasing as the viral potency $\xi$ decreases.  For example, close to the HB, which occurs at $\xi_{HB}\approx 0.042$, a partial eradication solution for $\xi=0.04$ gives a tumour burden $U\approx 2$. Note also the extension of the plateau of the periodic branch (green) close to the $U=100$ unstable equilibrium, before the branch point. This indicates that a ``square wave'' type of oscillations can be present for a moderately extended parameter interval in $\xi$. 

Although not shown in the diagram, the switch between the node and the spiral equilibrium typical of the partial eradication solution takes place along the red branch. For the chosen parameters in Fig.~\ref{Fbif}(a), this happens at $\hat{\xi} \approx 0.01675$, with spirals existing for a value $\xi$ such that $\hat{\xi}<\xi <\xi_{HB}$. Generally speaking and as shown in Fig.~\ref{nsregion}(b)-(c), the value at which the equilibrium type changes depends also on the other parameters $m,\gamma$ of the model. 

The system's behaviour also shows a strong, nonlinear dependence on viral death rate $\gamma$, as illustrated in Fig.~\ref{Fbif}(b). With respect to the case of $\xi$, the sensitivity of the model to $\gamma$ is somewhat reversed: intuitively, a surge in potency should act on the model in a similar way as a reduction in death rate and vice versa. For example, the branch of oscillatory solutions (green) out of the supercritical Hopf bifurcation (HB) shows an increasing maximum in $U$ as $\gamma$ decreases, opposite to what happens for $\xi$ (see the inset, in particular). 

The stable, impartial eradication solution branch (red) shows higher tumour volumes with increasing $\gamma$, and coalesces with the unstable $U=100$ branch (in black) at $\gamma = 1$.  For $\gamma>1$, the ineffective treatment solution is stable, as previously found from the analysis of the characteristic equation corresponding to this solution, i.e. Eq.~(\ref{ch-eq-partial}). A virus with a decay rate $\gamma > 1$ has no effect on the tumour. It is important to note that a mechanism identical to what we observe in the bifurcation plot for $\xi$ allows the existence of case (4) solutions, i.e. complete eradication. At a value of $\gamma \approx 0.0103$, the inset shows the termination of the oscillatory solutions (in green) and the occurence of a saddle-node point in the full eradication branch, making complete destruction of the tumour possible. From the biological perspective, this indicates that a right balance between the potency of the virus and its mortality must be achieved for an eradication to occur, depending on the growth rate $m$ of the tumour. In particular, as $\gamma$ is increased from zero, the model goes from full eradication to oscillations with an amplitude that decreases with $\gamma$, and subsequently to incomplete eradication up until $\gamma = 1$.

As previously mentioned, the full eradication solution gives rise to a singular Jacobian, making a purely numerical approach to continuation impossible. For solutions where $U\neq 0$, results have been obtained by using AUTO~\cite{AUTO07} and XPPAUT~\cite{Bard2002} softwares. For the case of solutions occurring for $U=0$, a combination of numerical methods and symmetry arguments have been employed. We assume that $U < I < V$, as exemplified by case (4) shown in Fig.~\ref{Fdyn}. If $\epsilon > 0$ and small, and we impose that $U\to \epsilon^n$, $V\to \epsilon^m$ and $I \to \epsilon^l$ with $n>m>l$, then the eigenvalues of the Jacobian $J$ can be numerically computed with an increasing approximation for growing $n,m$ and $l$. 

For example, in determining the stability of the full eradication branch in Fig.~\ref{Fbif}(a), we consider $U = 10^{-7}$, $V = 10^{-5}$ and $I = 10^{-4}$, substitute these values in the Jacobian and numerically evaluate the eigenvalues. For $\xi>\xi_{SN}\approx 0.0975$, all three eigenvalues turn out to be negative and real, whereas for $\xi<\xi_{SN}$ two are positive and one is negative. For example, choosing $\xi = 0.095$ gives eigenvalues $\lambda_1 \approx -0.15$, $\lambda_2 \approx -0.06$ and $\lambda_3 \approx 8\cdot 10^{-5}$. For the case $\xi = 0.099$, the first two eigenvalues are almost unchanged, but the last one changes sign and is $\lambda_3 \approx -2 \cdot 10^{-3}$. Similar results hold for the SN on the eradication branch for continuation in $\gamma$ (see Fig.~\ref{Fbif}(b)), and the method is consistent for all the parameters $m,\gamma$ and $\xi$ we have tested (not all shown here). These results have also been checked by integrating the equations of motion with XPPAUT, and confirming that the solution is indeed attracting when stable or repelling when unstable, as reported in the bifurcation diagrams.  

\begin{figure}[h!]
\centering
\includegraphics[width=0.49\textwidth]{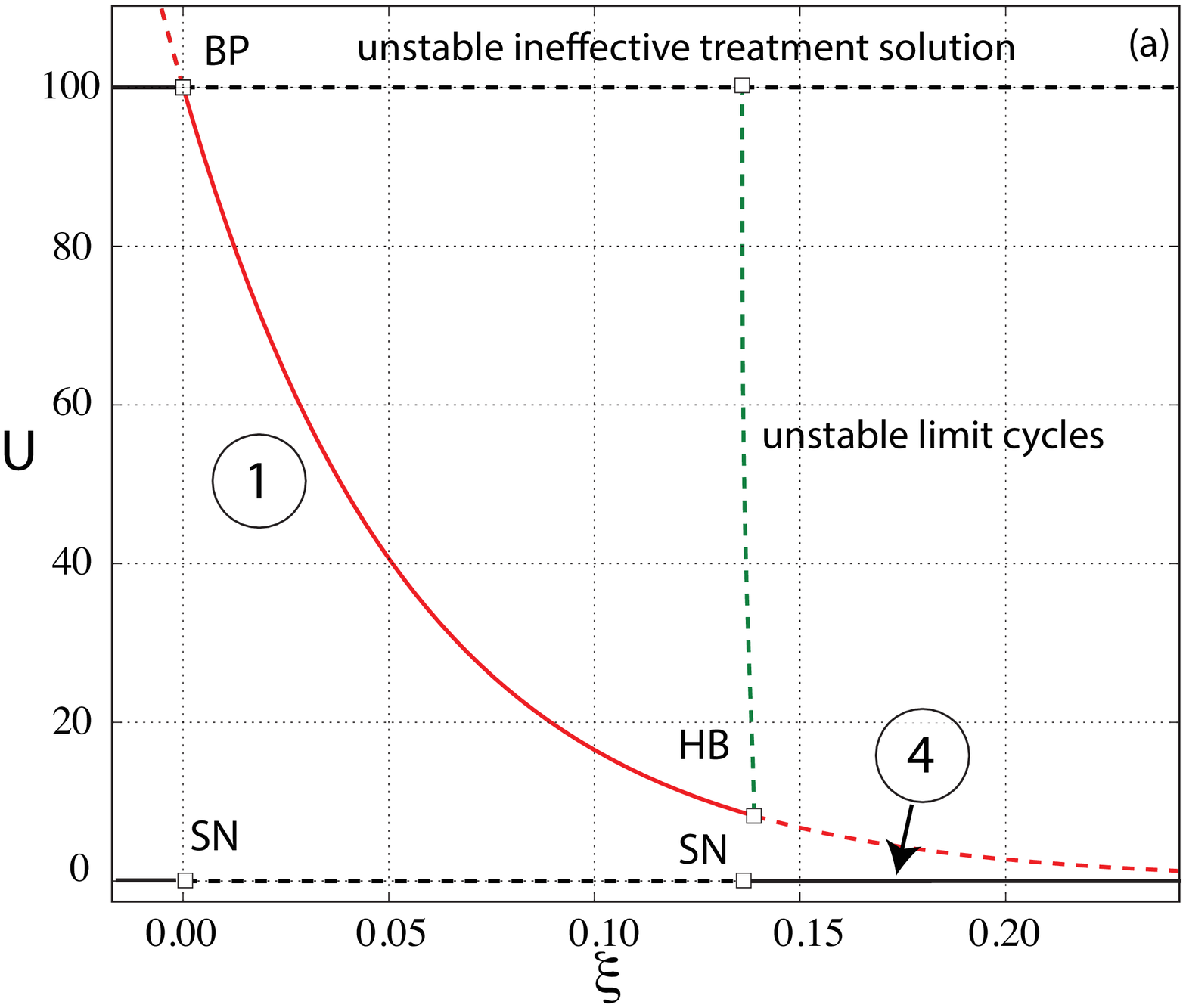}
\includegraphics[width=0.49\textwidth]{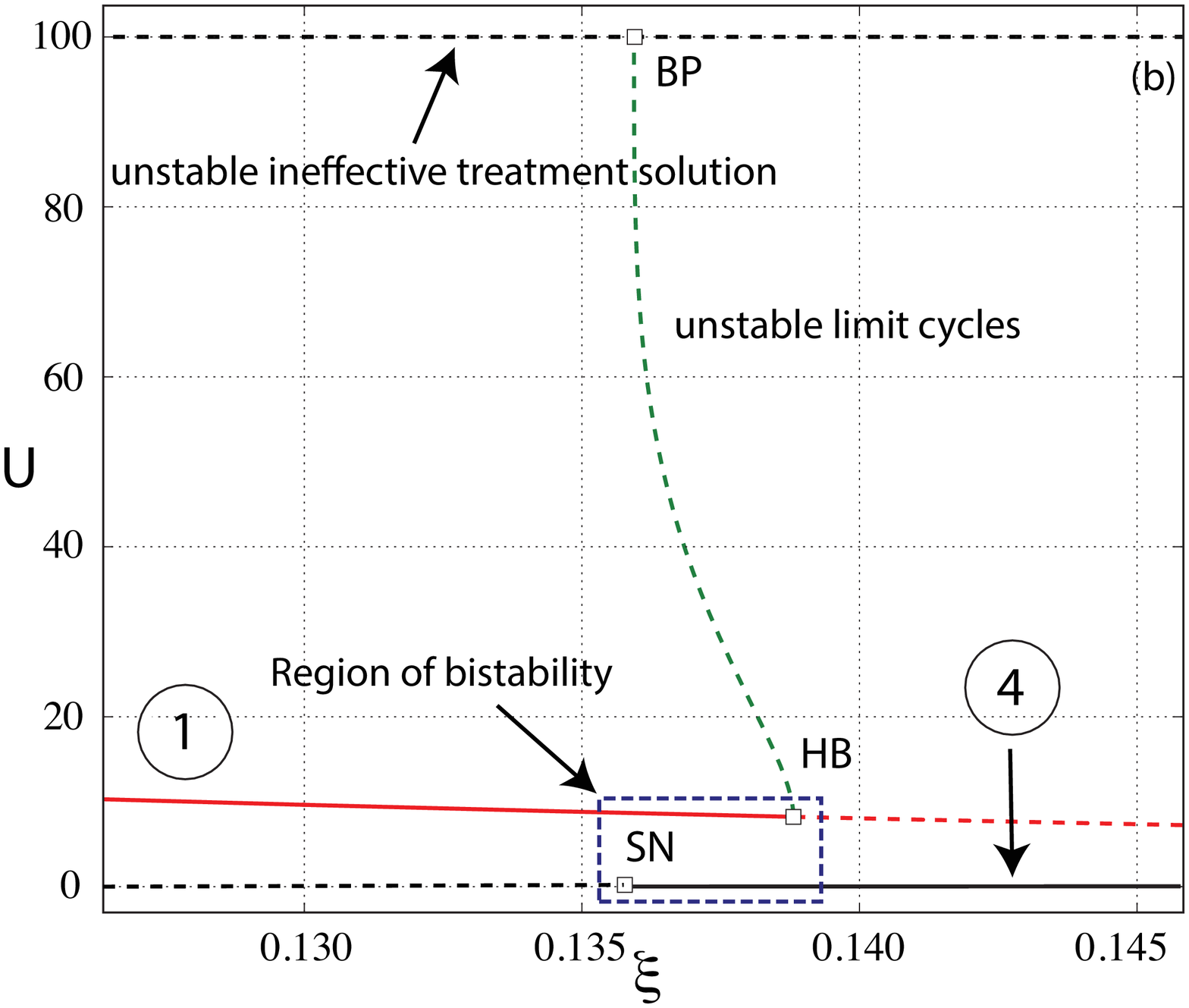}\\[5pt]
\includegraphics[width=0.49\textwidth]{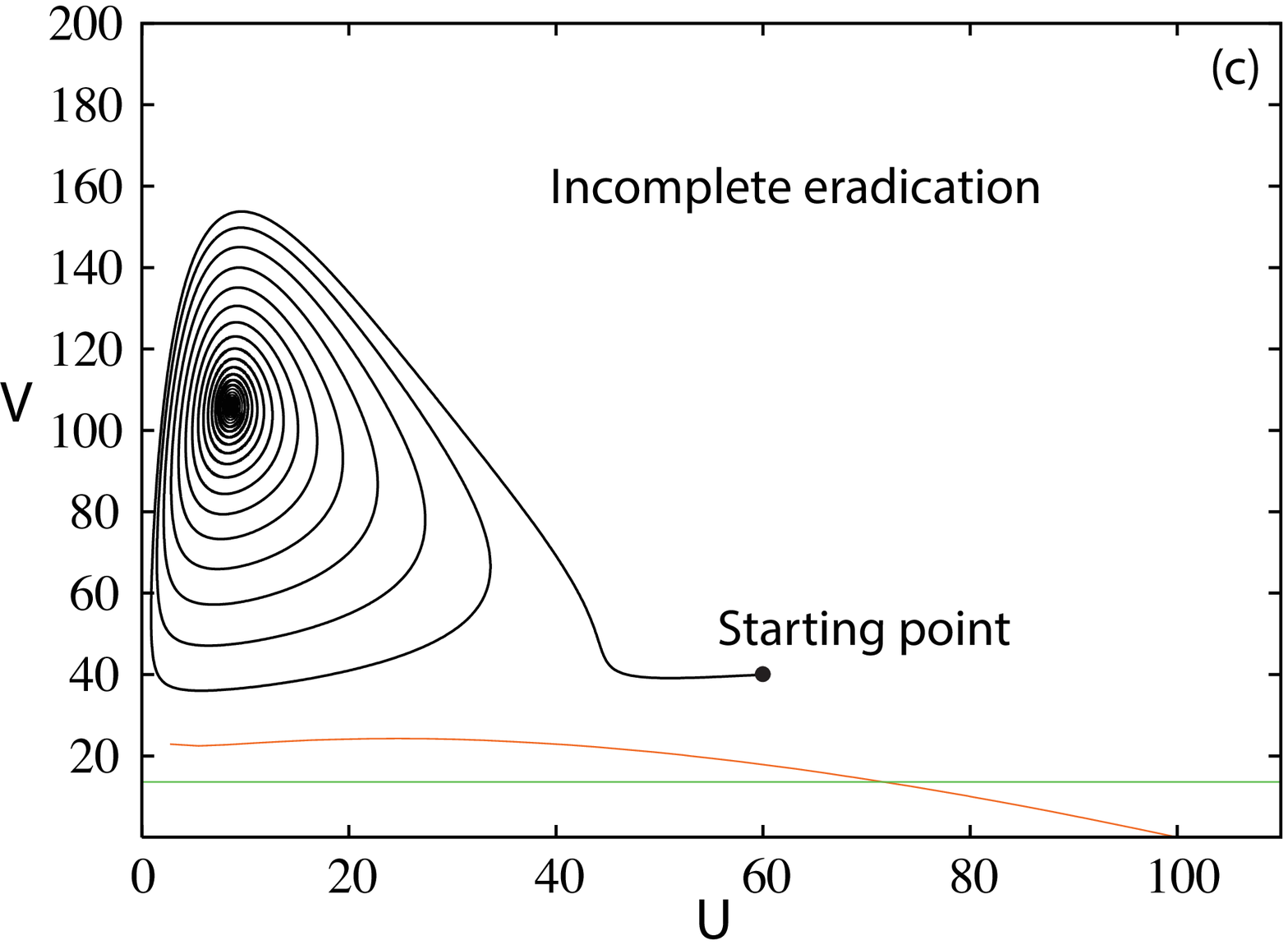}
\includegraphics[width=0.49\textwidth]{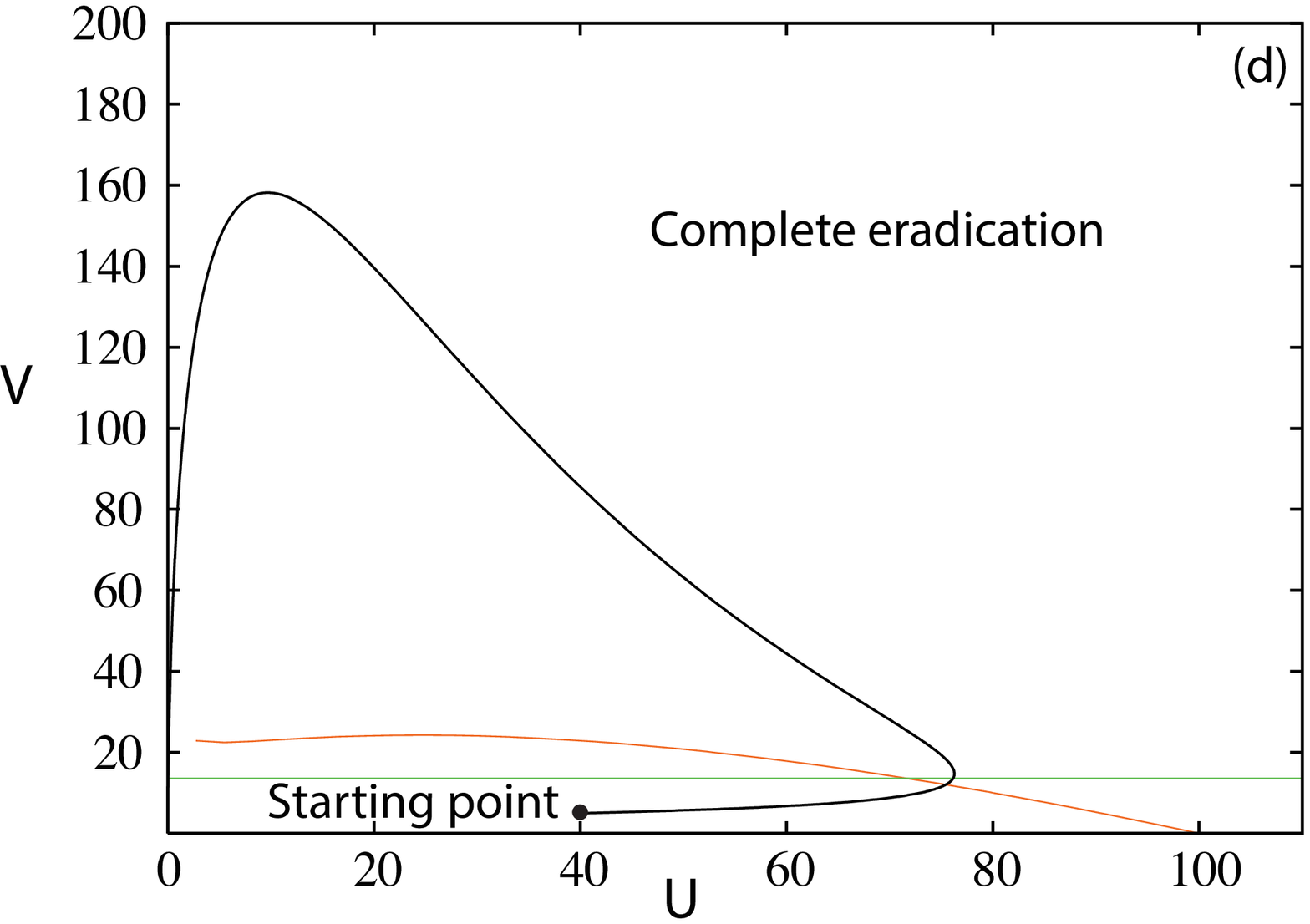}
\caption{Bifurcation plots and bistable solutions for fixed parameter values $m=0.5$, $\gamma=0.1$. The rectangle in (b) shows the area where two solutions of different nature coexist, delimited by $\xi_{SN} \approx 0.1359$ and $\xi_{HB} \approx 0.1388$. A spiralling solution to an incomplete eradication is shown in (c) and occurs for initial conditions $U(0)=60$, $I(0) = 10$, $V(0)=40$, for a parameter $\xi_{SN} < \xi = 0.136 < \xi_{HB}$ . A fully eradicated solution is shown in (d) and instead occurs for $U(0)=40$, $I(0) = 10$, $V(0)=5$, for the same value $\xi = 0.136$. Nullclines, i.e. the loci of points corresponding to $U'=0$ and $V'=0$, are in red and green.}\label{Fbista}
\end{figure}

One important feature of the model is that it does not support stable oscillations for all biologically meaningful combinations of parameters. For some choices, a different structure of bifurcation plots emerges, with significant consequences from the biological perspective. In this sense, a typical example for a continuation in $\xi$ is illustrated in Fig.~\ref{Fbista}(a). An unstable periodic branch (green)  originates from a subcritical Hopf bifurcation (HB) and maintains its unstable character until it collapses with the $U=K=100$ (black) branch. For this diagram, viral potency $\gamma$ is the same as in Fig.~\ref{Fbif}(b), but a value of $m=0.5$ (moderately high growth rate) is chosen, whereas both previous diagrams have been obtained with a $m=0.1$ (moderate growth rate). A more aggressive tumour, assuming that the potency of the virus is the same, does not engage in oscillatory behaviour with the virus, but only partial or full eradication are possible (i.e. black and red lines). 

It is interesting to stress that in this case, as shown in Fig.~\ref{Fbista}(b), the saddle-node (SN) on the full eradication $U=0$ branch (in black) occurs for a value $\xi_{SN}$ that is less than the value $\xi_{HB}$ at which the subcritical Hopf (HB) originates. This occurrence is due to the fact that the periodic branch shows increasing values of max $U$ for decreasing values of $\xi$ when it is unstable. This is the opposite of what happens for the stable periodic branch described in Fig.~\ref{Fbif}(a), where $\xi_{HB} < \xi_{SN}$ and the stability of the eradicated solution does not switch in this way. 

The change in the order in which the SN and HB emerge as $\xi$ is increased is responsible for the generation of a region of bistability, where two separate and distinct equilibria exist for an interval of potency values. For values of $\xi$ in this region, different initial conditions can lead to different outcomes, as shown in Fig.~\ref{Fbista}(c)-(d): the initial dosage of viral load and the numbers of infected and uninfected tumour cells can strongly influence the final fate of the system and, as it will be clear shortly, lead to somewhat unexpected results. In the first case (Fig.~\ref{Fbista}(c)), a spiralling solution achieves an incomplete eradication, which belongs to the red branch in Fig.~\ref{Fbista}(b). Conversely, the second case shows a complete eradication to a vanishing tumour, after traversing two maxima in $U$ and $V$ respectively, corresponding to the black branch in Fig.~\ref{Fbista}(b). A small variation in the initial conditions can hence result in the therapy being effective or instead giving rise to a partial eradication.

\begin{figure}[h!]
\centering
\includegraphics[width=.45\textwidth]{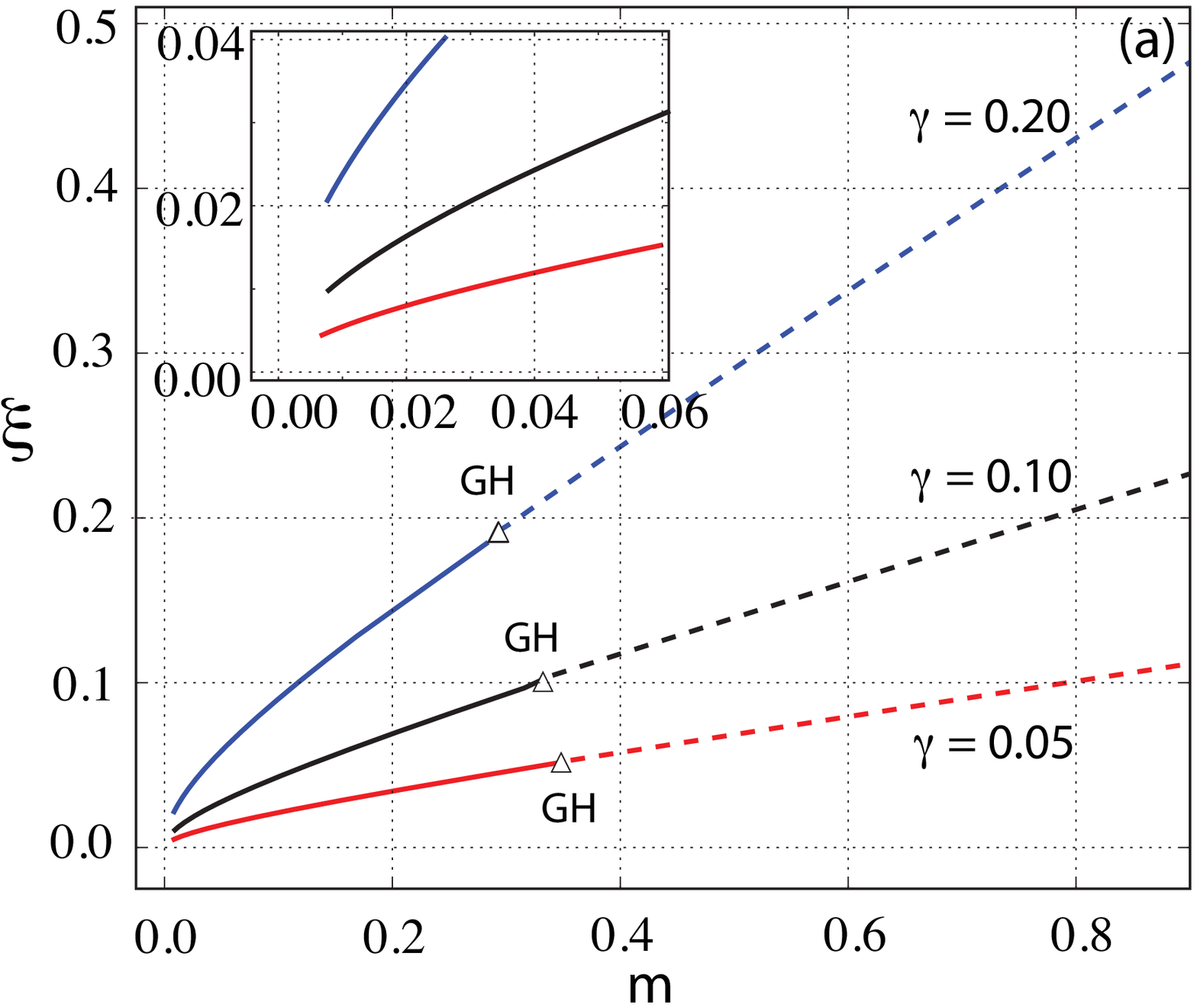}
\includegraphics[width=.45\textwidth]{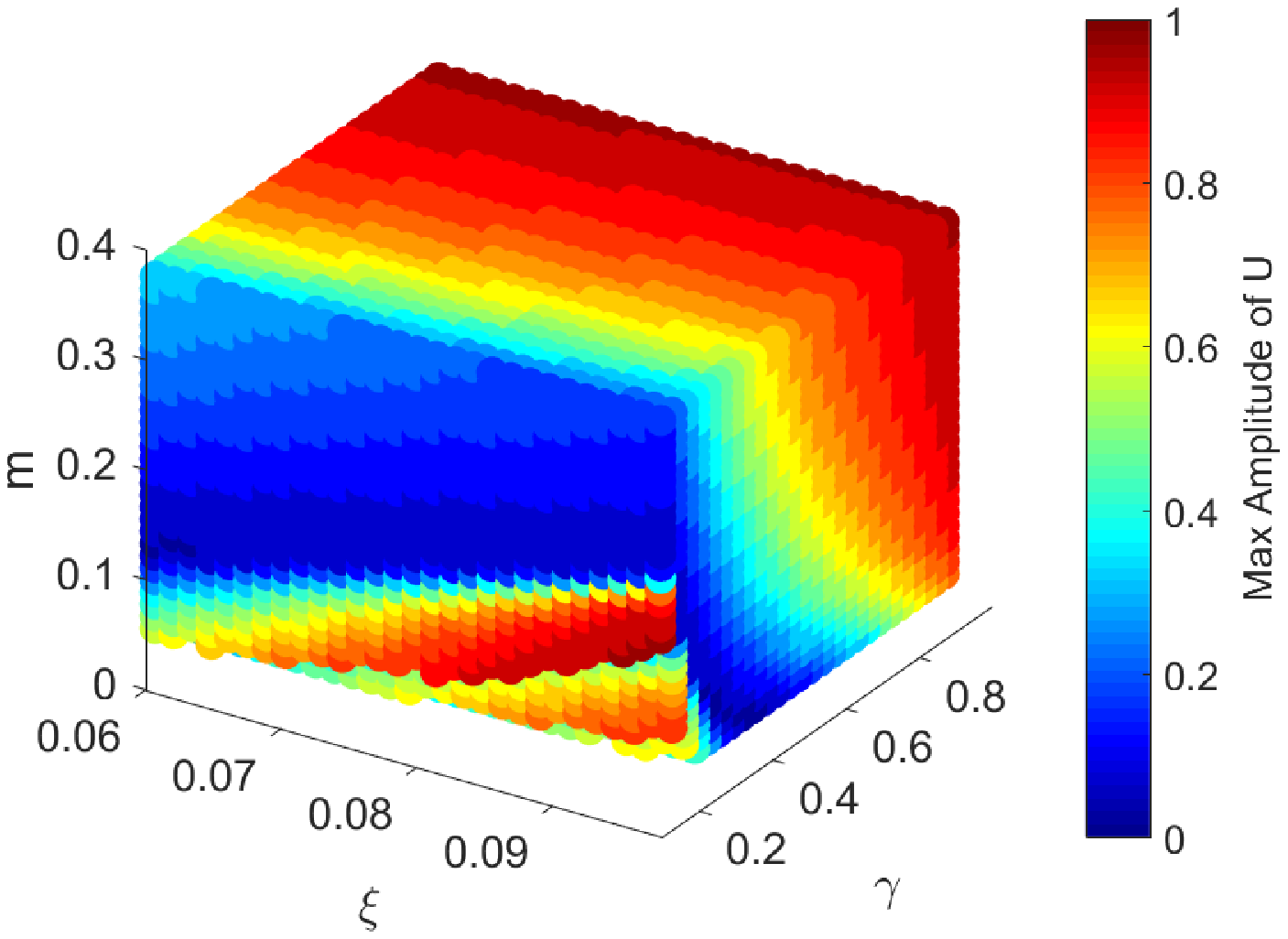}
\caption{Different, two-parameters continuations in (a) for $m$ and $\xi$ for branches of Hopf bifurcations at different values of $\gamma$. Branches of supercritical Hopf bifurcations are shown in continuous lines, whereas those for subcritical bifurcations are in dashed lines. Generalised Hopf points are indicated by GH. Note that the branches cease to exist for low values of $(m, \xi)$,  indicating the system cannot support either stable or unstable oscillations when parameters are sufficiently small (see the inset). Plot of the corresponding amplitude of stable limit cycles for points in the $\xi,m,\gamma$ parameter space are in (b). The colour of the point corresponds to the maximal value of the amplitude of the limit cycle in $U$.}\label{2par}
\end{figure}

The existence of this area of bistability is associated with the presence of a subcritical Hopf bifurcation whose loci of points in $\xi$ and $m$, and for different values of $\gamma$, are plotted in Fig.~\ref{2par}(a). Generalised Hopf points (GH) separate subcritical Hopf points (dashed lines) from supercritical Hopf bifurcations (continuous lines). Note that, if the the growth $m$ is sufficiently small, no Hopf bifurcation can be present and the system does not support oscillations, either stable or unstable. For example and as a result of the interruption of the Hopf branches shown in the inset of Fig.~\ref{2par}(a), any one-parameter bifurcation plot in $\xi$ for a fixed $\gamma =0.1$ and values of $m \lessapprox 0.008$, does not contain a stable or unstable oscillatory branch, since no Hopf point exists for such values. Biologically this indicates that there is a lower bound on the tumour growth rate for oscillations (stable or unstable) to exist, implying that a very slow growth in general leads to a complete eradication for sufficiently potent virus, as long as its death rate is not excessively pronounced. 

A numerical analysis of the model for a range of $\xi,\gamma$ and $m$ values shows that limit cycle amplitudes for $U$ do not follow a clear pattern, as captured by Fig.~\ref{2par}(b). Oscillations of different amplitudes can be achieved by the system and, as observed previously, depending on the growth rate of the tumour they can be enhanced by increases in $\xi$ and decreases in $\gamma$. 

\section{The effect of dosage applications and their optimisation}

As shown, there are large sections of parameter space that give birth to regimes with dormant tumours or tumour-virus oscillations, which can give rise to different outcomes when coupled with clinical therapies. For the purpose of the present study, one of the simplest possible therapeutic practices is considered: the administration of constant dosages of viral loads via external injections and at given time intervals. If the treatment is over the course of $n$ injections with $\kappa$ number of days between injections, a virus injection protocol $u_V(t)$ can be summarised by the following, generic schedule: 
\begin{equation}
u_V(t) = \left\{\begin{array}{cc} \displaystyle\frac{D_0}{n} \qquad & \qquad t = (i-1)\kappa,  \hspace{0.5cm} 
\text{where }i = 1, \ldots, n, \\[6pt] 0\qquad & \text{otherwise.} \end{array} \right.\label{Eq5}
\end{equation}
\noindent 
Given this simple scheme, let us consider how dosage perturbations affect regions of the bifurcation diagrams and if they can result in either tumour eradication or a stable tumour size below a given threshold. The two typical scenarios we consider are oscillations and bistability.

\subsection{Effects of injections on a stable, oscillatory trajectory}

After exploring different areas of the parameter space that give rise to oscillations, the main finding reported in this study is that simple therapies given by Eq.~(\ref{Eq5}) do not alter the long term behaviour of the model independently of the amplitude or period of oscillations. If an oscillatory, stable state exists between virus and tumours, an increment in the viral load through injections does not achieve complete eradication to zero. From the dynamical point of view, an increase in viral load via external perturbation cannot force the system out of the basin of attraction of a stable limit cycle. Nonetheless, transient phenomena do exist and are worth discussing.

\begin{figure}[h!]
\includegraphics[scale=.45]{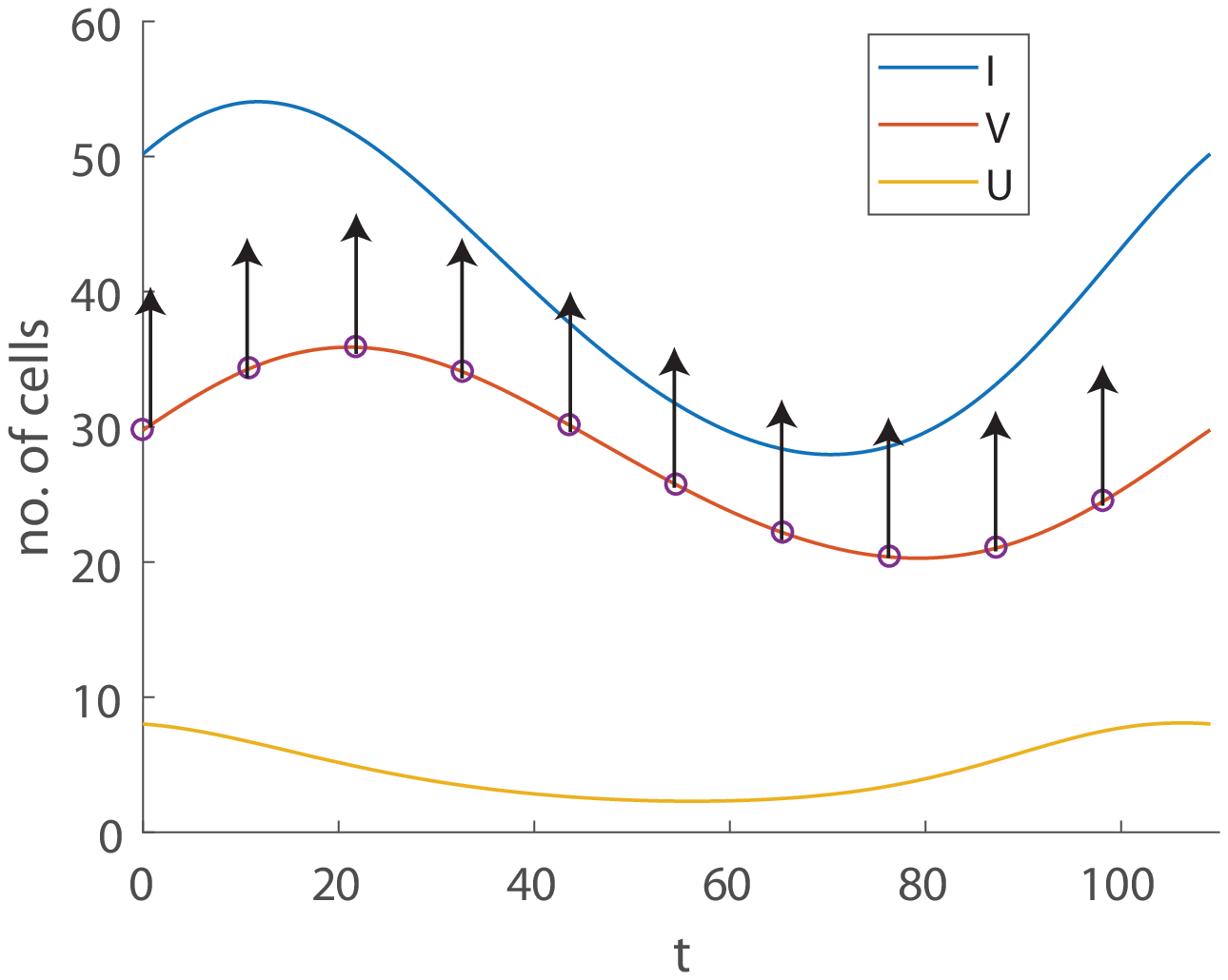}
\includegraphics[scale=.45]{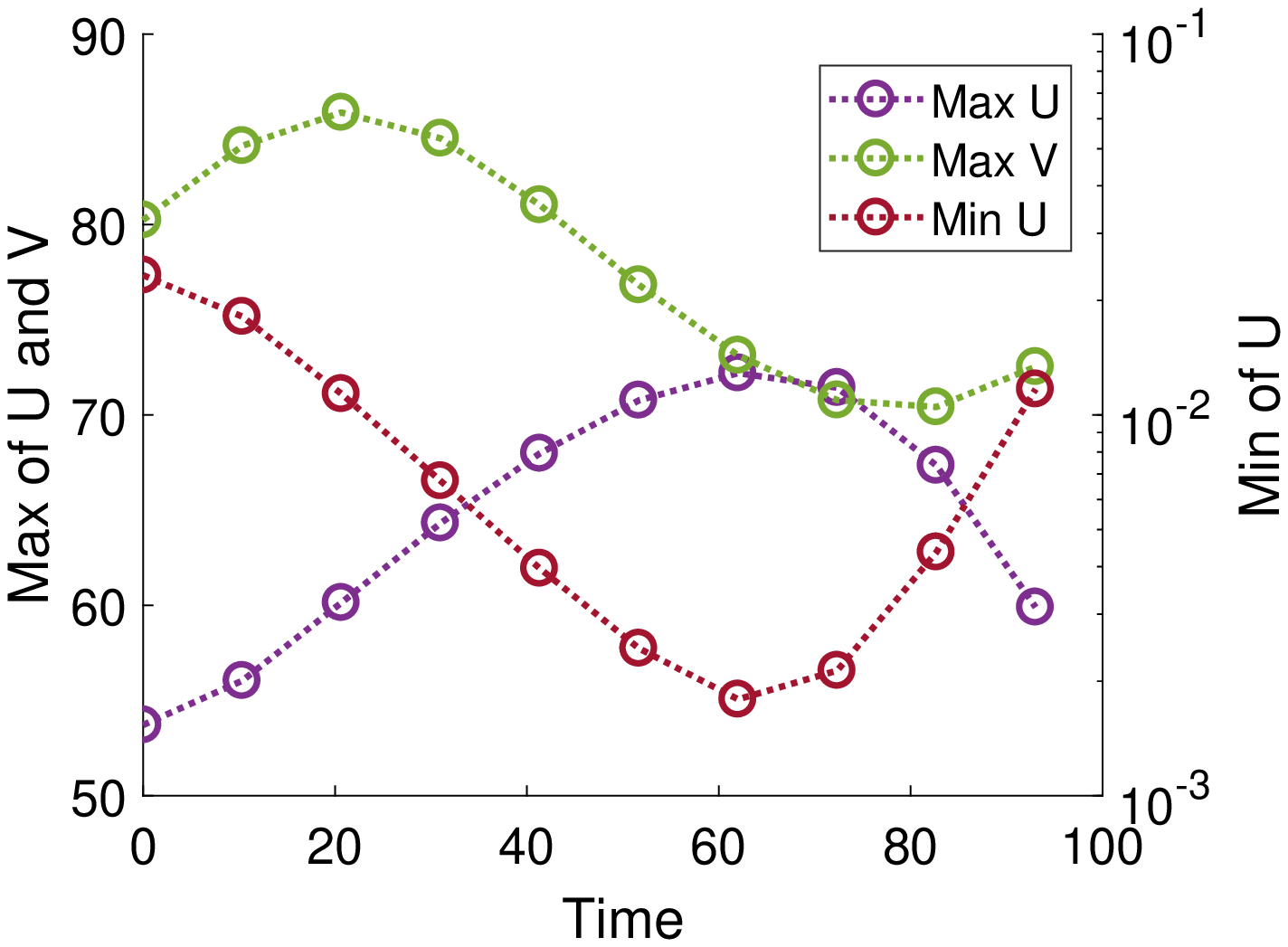}\\[5pt]
\includegraphics[scale=.45]{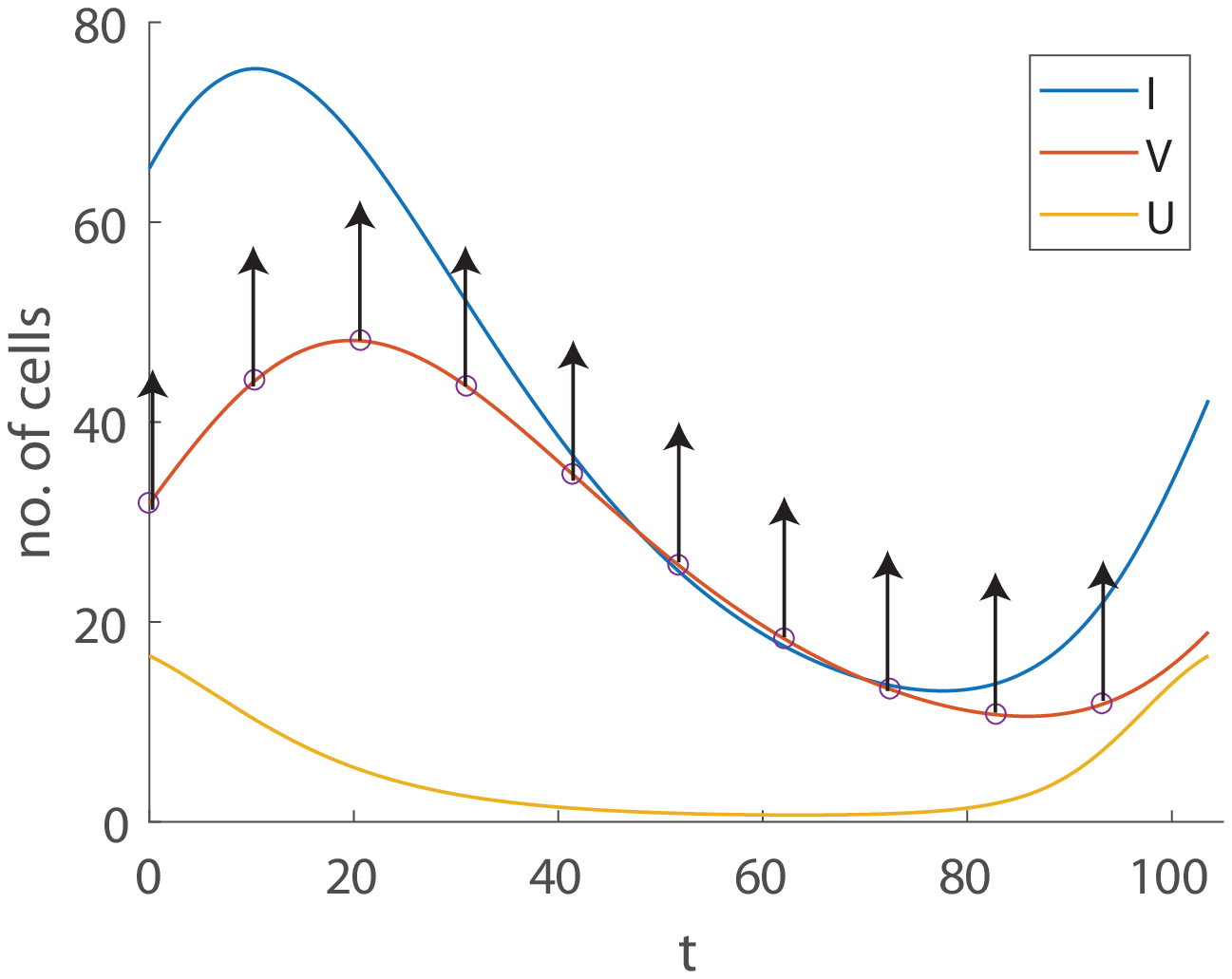}
\includegraphics[scale=.45]{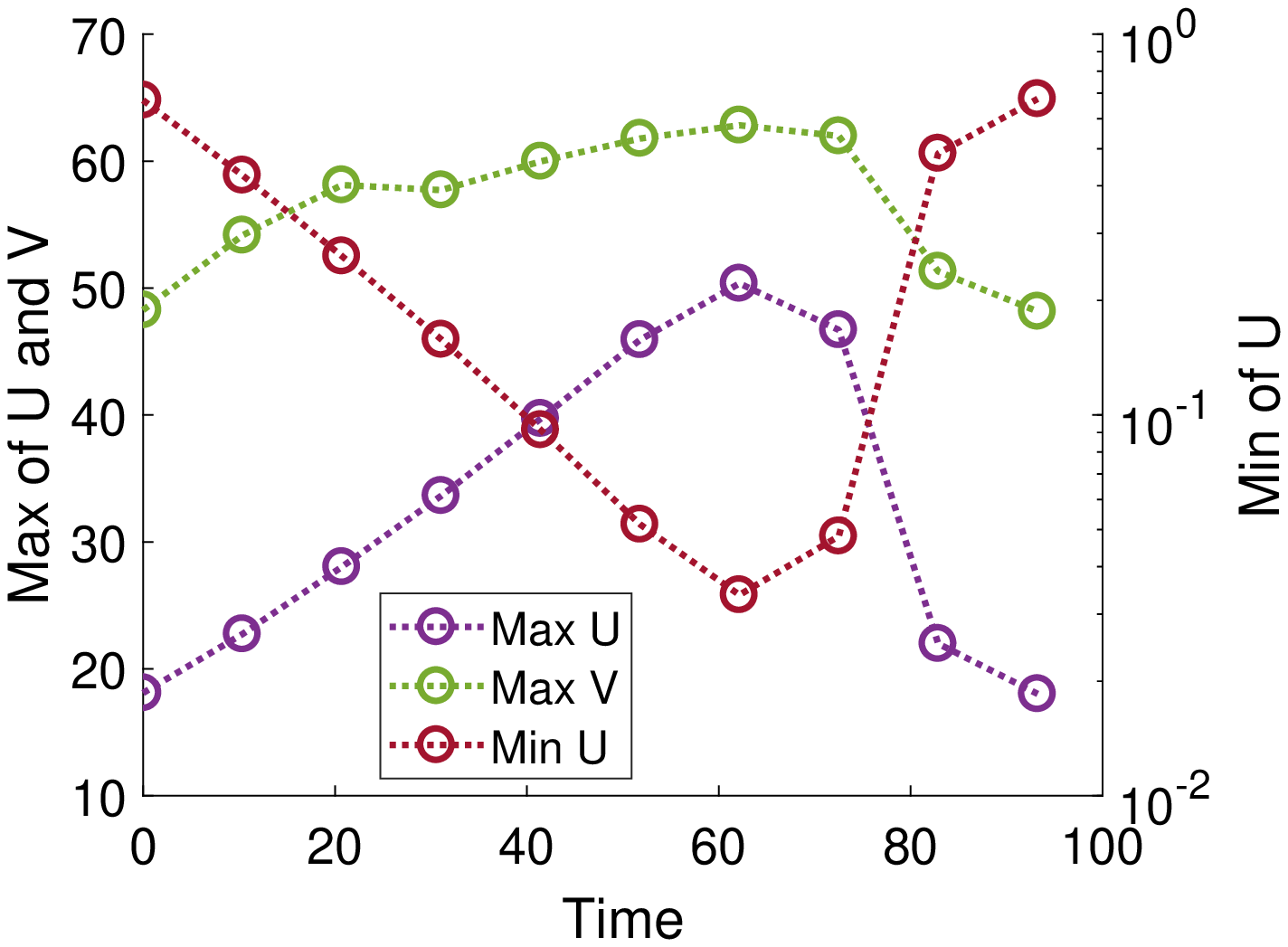}
\caption{Perturbations in the days between two treatments $\kappa$. Two different limit cycle regimes have been plotted for $\gamma = 0.1, m = 0.2$ (a) $\xi = 0.06915$ and (c) $\xi = 0.06993$. The maximum and minimum uninfected cell number is plotted in (b) and (d) for the corresponding value of $\kappa$ represented by an upward arrow in (a) and (c). Note the different scales, since the oscillations considered in (a) and (c) have different amplitudes.}
\label{D0sims}
\end{figure}

Let us consider two injections, i.e. $n=2$, for a system already in a stable oscillatory state. The number of days $\kappa$ between injections alters the size of the tumour and virus populations as the system returns to its stable state. In Fig.~\ref{D0sims}, a single period of two different stable limit cycles of the model is shown, with arrows representing the instants at which injections that increase the viral load have been administered. The corresponding maximum and minimum tumour size, along with the maximum virus count reached, is also presented. 

Injections that occur at different phases of the cycles have different outcomes. In particular and for large or small oscillations, as can be seen in Fig.~\ref{D0sims}(b)-(d) for the red and magenta curves around $t\approx 62$, dosing the virus close to the minimum in tumour population provides a typical outcome: the tumour initially responds to the injection by undergoing the lowest resulting minimum, but this is followed by a rebounds that causes $U$ to reach the highest value (max $U$ in the plot) of all other tested injections. Note that, in some cases and for sufficiently high dosages, the minima achieved by $U$ can be pushed to values so low to become experimentally undetectable. Injections at other instants within one oscillation period yield rebounds proportional to the original amplitude of the limit cycle, with best results occurring for the lowest amplitudes.

Perturbing the number of days between injections $\kappa$, the total injection amount $D_0$ or the number of injections $n$ does not affect the long term dynamics (not shown here), which remains oscillatory in the long term.

\subsection{Effects of injections on a trajectory in the bistable region}

For a solution in a bistable region, the final outcome of any injection is highly dependent on the initial tumour size and viral load. In particular, due to the complex structure of the basin of attraction of the two competing solutions, i.e. full eradication and an incomplete quiescent state, doses that are higher than a specific threshold, which is in turn highly dependent on the system parameters, can lead to a partial eradication rather than a complete one.

\begin{figure}[h!]
\centering
\includegraphics[width=0.48\textwidth]{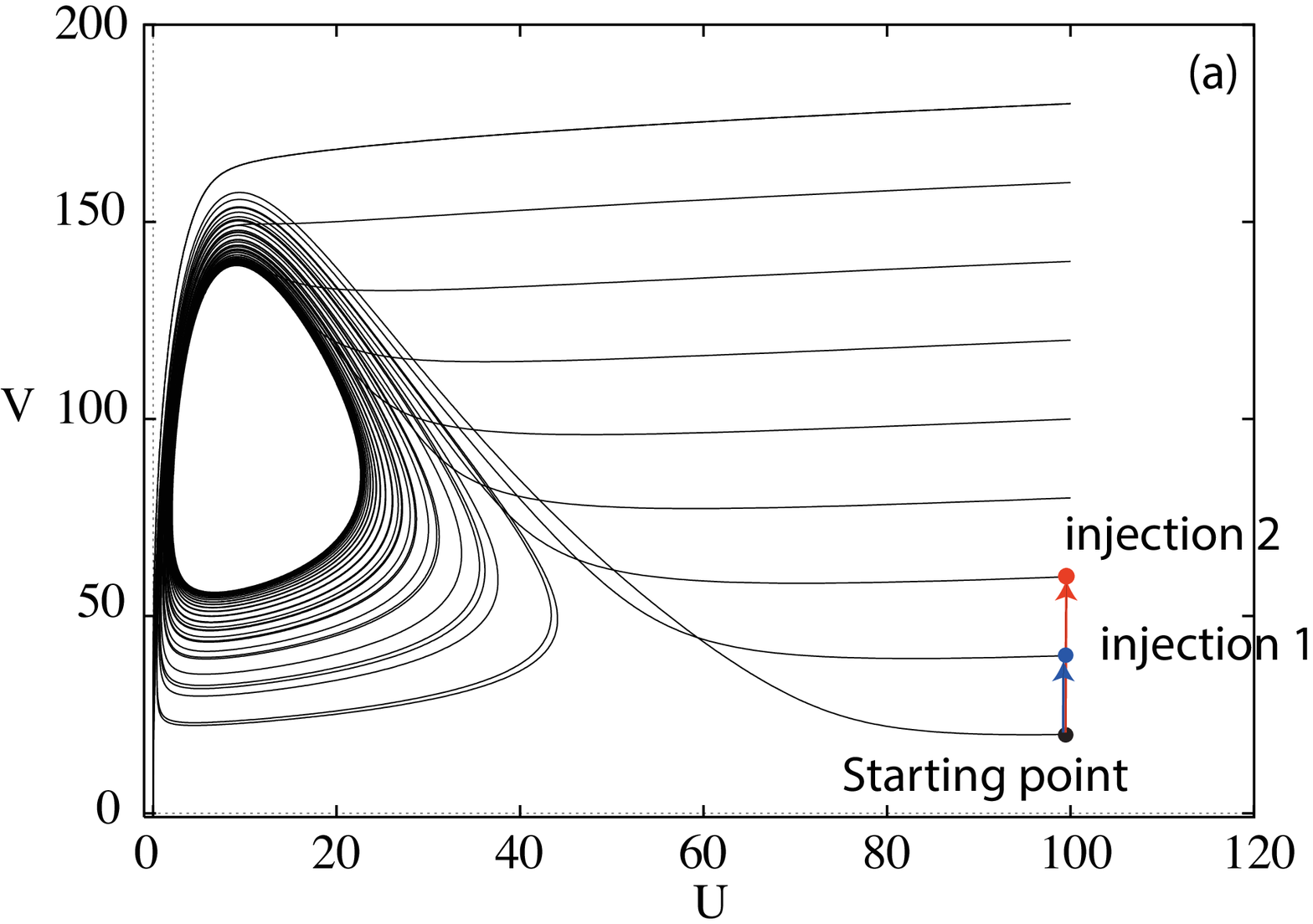}
\includegraphics[width=0.48\textwidth]{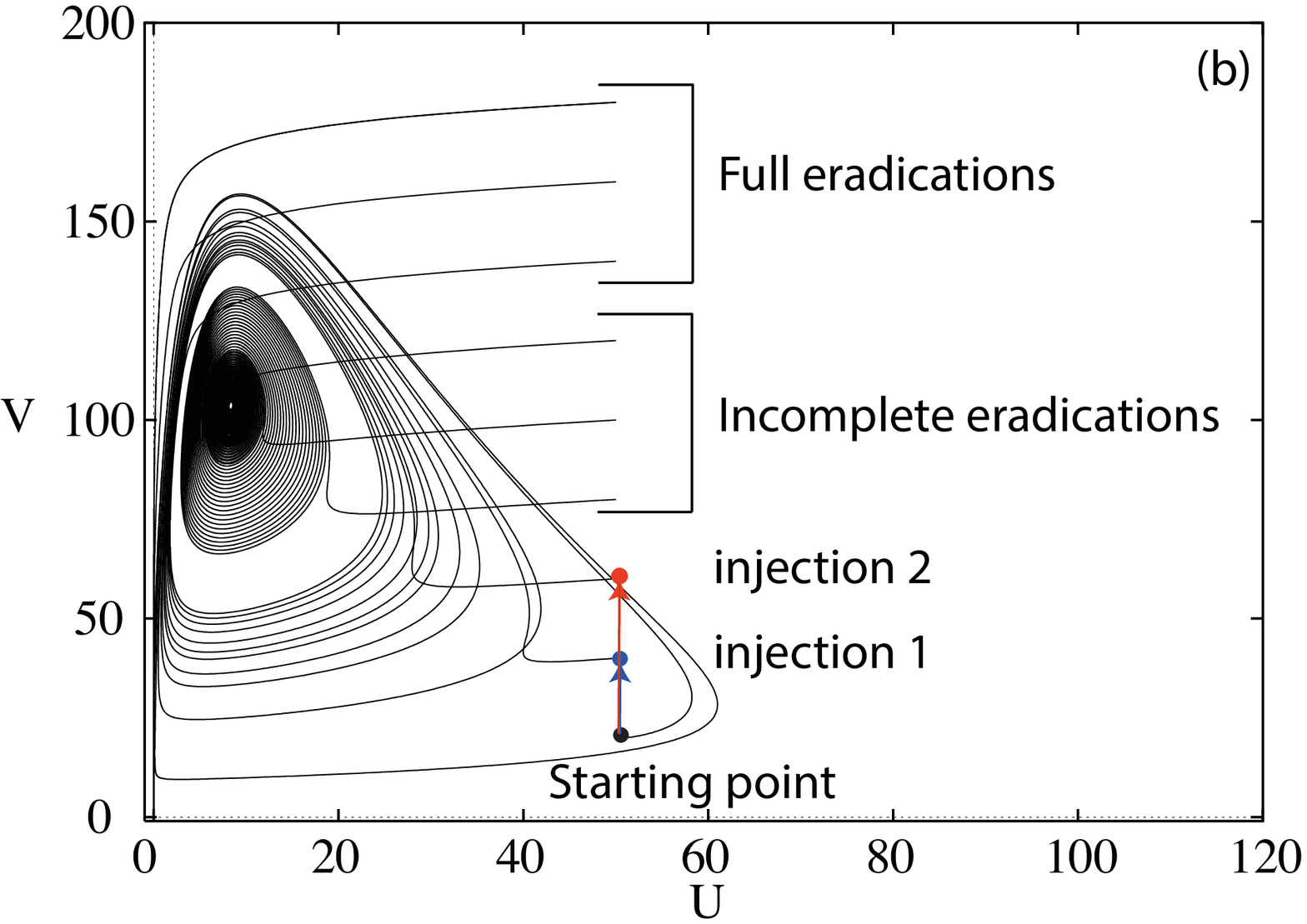}
\caption{Typical cases of dependence on injected viral dosage $D_0$ for a system in a bistable scenario. Examples of two injections with increasing dosage (i.e. injections 1 and 2) are also sketched. The effect of these injections is to push the starting point to larger values of $V(0)$, depending on the dose that is administered. For the same initial tumour size, different dosages result in either tumour eradication or tumour stabilisation. Initial fixed conditions in (a) are given by $U(0)=100$, $I(0)=10$ and by $U(0)=50$, $I(0)=10$ for (b). In both cases, $V(0)$ varies from a minimum of $20$ to a maximum of $120$ in constant steps and the parameters are $m = 0.5$, $\gamma =0.1$ and $\xi = 0.138$.} \label{bistabsim}
\end{figure}

As typical scenarios, consider the administration of single injections of increasing dosage as depicted in Fig.~\ref{bistabsim}. Depending on the initial uninfected tumour $U(0)$, injections can lead to different outcomes or even have no effect on the final state. Considering the case of a high tumour size (Fig.~\ref{bistabsim}(a), $U(0) = 100$), different dosages always result in final eradication. Some dosages can lead to transient oscillations in the $U-V$ plane, but eventually eradication is achieved for all plotted trajectories. If, instead, the initial tumour size is smaller (Fig.~\ref{bistabsim}(b), $U(0) = 50$), a full eradication is obtained only if the dose is either sufficiently low or sufficiently high, whereas there is a considerable interval of possible doses that push the system to a stable spiral corresponding to a dormant state, where eradication is not complete. Note that the first two low dosage injections, i.e. injections 1 and 2 in Fig.~\ref{bistabsim}(b), also lead to a final eradication state after few oscillations on the $U-V$ plane. 

This result is interesting, as it suggests that, for given initial tumour size and characteristics of the virus, there can be a unique interval of dosage sizes that does not result in treatment success. Boosting the amount of virus does not always guarantee a successful outcome.
 
\section{Discussion}

The model proposed in this work shows a number of interesting features, both from the mathematical and the biological points of view. Firstly, a range of possible dynamical outcomes, based on the value of the model parameters and, in some cases, of the initial conditions, have been found. A number of nontrivial bifurcation scenarios have also emerged, with the presence of an important system equilibrium (i.e. full tumour eradication) that is characterised by a singular Jacobian. This occurrence has required the use of a hybrid combination of numerical continuation, symmetry considerations and integration of the model to map out the dynamics as a function of relevant model parameters.

The model provides a few insights into the interactions between an oncolytic virus and a tumour growing with a realistic proliferation law. One of the main limitations of the present approach is the endless influx of viral load that occurs in the model: once the viral cycle is set into motion, and unless viral death rate is excessive (i.e. $\gamma > 1$), there is no natural stopping mechanism for viral infections to continue endlessly. This simplification is, for example, responsible for the appearance of dormant, partially eradicated tumours, which, after an initial transient, perpetually coexist with a constant viral load. These dynamics are common for models with unlimited reservoirs of populations~\cite{Wilkie2013201} .

Another important constraint is represented by the limited number of parameters used and their inherent inability to fully account for tumour-virus dynamics in detail. For example, we have introduced two general terms $\xi$ and $\gamma$ that aim to capture the potency and death rate of the virus and depend on the virus infectivity rate, $\beta$. These parameters are meant to encapsulate a large variety of different viral characteristics and can be associated to features as diverse as burst size, reproduction rate, spreading ability, and diffusivity. A similar observation must be made for the growth rate $m$: this value condenses a large number of often independent and highly variable features of tumour growth, which are highly sensitive to nutrients, vascularisation, extra cellular matrix characteristics, and so on.  

Notwithstanding these limits, the model shows that, for a given rate of growth, a tumour responds in different ways to viral particles that have different, generic characteristics. As shown in Fig.~\ref{Fbif}, an increase in viral potency $\xi$ or a decrease in viral death rate $\gamma$ drives the system through similar stages of typical dynamics, from partial eradication to tumour-virus oscillations. At sufficiently large values of $\xi$ and $\gamma$, for instance $\gamma > 1$, the scenarios are instead opposite, with full eradication and inefficient treatment, respectively. 

A metastable regime that appears somewhat counterintuitive is represented by the so-called ``square-wave'' oscillations, which are observed in a small interval of biologically relevant parameters (see Fig.~\ref{Fdyn}(c)). Given the limitations of the model proposed here and the size of the parameter space where this dynamics takes place, it may be unlikely that such an extreme tumour expansions can be directly observed in a clinical setting. Nonetheless, the switch between a quasi-eradicated to a quasi-ineffective treatment regime points to the importance of achieving a complete wipe out of the tumour if a sudden resurgence is to be avoided. 

The existence of an extended area of the parameter space where oscillations among system variables arise is also worth noticing. These regimes, which also tend to respond nonlinearly to external injections (see Fig.~\ref{bistabsim}), have been known for quite some time in clinical settings {\bf XXX add REF if existing - Adri, please}. One major finding for this model is that virotherapy can prevent oscillations from occurring if the potency is sufficiently strong or, alternatively, the virus tends to survive for sufficiently long times in the infected population. Furthermore, and this is particularly interesting, oscillations tend to have larger amplitudes and periods for increasing $\xi$ (or decreasing $\gamma$), before they disappear completely for sufficiently high (or low) values. This is worth reflecting on, especially from the clinical perspective. Designing a potent virus that is still not sufficiently resilient may turn out to be a riskier strategy, since it could trigger larger fluctuations in the tumour population. These oscillations also occur at relatively distant time intervals from each other and long periods of tumour inactivity may be misinterpreted as successful eradication. Looking at Fig.~\ref{Fbif}(a) and assuming that a low value of uninfected tumour cells $U$ represents a good outcome, a less potent virus, say with a $\xi \approx 0.04$, results in a quiescent tumour of a smaller size than the amplitude of oscillations caused by a highly potent virus with, for instance $\xi \approx 0.08$ (i.e. twice as potent). This is also true from the point of view of resilience, see in particular the inset of Fig.~\ref{Fbif}(b): a virus that remains active for longer, say $\gamma \approx 0.015$, produces oscillations with very high values of $U$, whereas a virus decaying twice as fast, say with $\gamma \approx 0.03$, produces a stable, silent tumour of a smaller size. All this shows that therapeutic strategies must be chosen carefully and thoughtfully, and that optimal design of an oncolytic virus must be targeted on the tumour characteristics, in particularly its proliferation rate. It could be quite interesting to test these theoretical findings in vitro. 

Note also that, even when external interventions with extra viral dosages are taken into account, the answers provided by our analysis do not appear trivial. Firstly, the existence of bistability and dependence from initial conditions has important effects. As seen in Fig.~\ref{bistabsim}, different initial viral loads and dosages can result in different outcomes, often in unpredictable ways. It is not true that a larger initial viral load always results in eradication: there is a large interval of values of dosages for which eradication is not possible and, quite interestingly, the system privileges either a sufficiently high or sufficiently low viral load for eradication. Starting at a smaller viral load is successful because it first allows the tumour to initially grow to a larger size, which thus elicits a stronger viral response. This response can wipe out the tumour completely, with no risk of ending in a dormant phase. Although this feature has been observed previously, for example, in simpler systems in tumour-immune dynamics and predator-prey models~\cite{Davis1962, frascoli2013}, it is the first time, as far as the authors are aware, that it is noted in virotherapy modelling. Clearly, the fact that our model hypothesises that the virus can penetrate the tumour and diffuse within its cells with no hindrances, has to be taken into consideration and is one of the drivers of this effect. Notwithstanding this, the result points to the existence of a preferred threshold in the size, for some values of the system parameters, where a limited quantity of viral load is preferable over a larger amount. Although strategies resulting even in a controlled and partial growth of a tumour have to be evaluated and considered with extreme care, the fact that a low viral load can still produce positive outcomes should be investigated further in laboratory and clinical settings. We remind the reader that the present model does not allow for a thorough description of the dynamics of virus penetration and diffusion, which certainly play a fundamental role in the success of virotherapy.

Secondly, therapies that couple with external injections of viral loads could have very different outcomes depending on the state of the system. Not only, as we have just highlighted, they can perturb a trajectory that was meant to be of full eradication into a dormant state, but, as shown for oscillations in Fig.~\ref{D0sims}, they can have a transient, often negative effect on the whole system. If administered when the system resides on a stable oscillating state, these injections, depending on when in the cycle are provided, tend to increase the amplitude of few cycles of oscillations before the system goes back to its original fluctuations, with no ability of driving the model out of this phase. There is generally no positive relevant effect in reducing the magnitude of periodic behaviour in the long term. Strategies that instead optimise the quality of the oncolytic virus seem to be preferable, as (see Fig.~\ref{Fbif}) oscillations can be reduced or damped to zero either by increasing the potency or the life span of virus at the right amount.

In this sense and as also Fig.~\ref{2par}(a) explains, the finding that oscillations that exist for different values of the parameters are suppressed when the growth rate $m$ is sufficiently small is very relevant and informative for therapeutic choices. Rather than complex injection schedules or larger amounts of externally provided virus, this model seems to promote pharmacological interventions that aim at blocking or reducing the growth of the tumour. It will also be interesting, for further studies, to establish whether combination therapies or interventions specifically targeted at boosting the immune response (not modelled here) could also improve outcomes, and how changes in the diffusion and penetration efficiency of infection waves may change the trends observed in the present model.

\section*{Acknowledgements}

ALJ, FF and PSK gratefully acknowledge support for this work through the Australian Research Council Discovery Project DP180101512, ``Dynamical systems theory and mathematical modelling of viral infections''.

\appendix
\section{Non-dimensionalisation of model}
\label{sec:appendix}

Here we discuss the non-dimensionalisation of Eqs.~(\ref{AEqs1})-(\ref{AEqs4}) to obtain the system in Eqs.~(\ref{E4})-(\ref{E6}). First consider the original model:

\begin{align*} 
\frac{dU}{d\tau} & = r\ln\left(\frac{K}{U}\right)U -\frac{\beta U\hat{V}}{U+I},\\
\frac{dI}{d\tau} & = \frac{\beta U\hat{V}}{U+I}-d_II,\\
\frac{d\hat{V}}{d\tau} & = -d_V\hat{V}+\alpha d_I I. 
\end{align*} 

\noindent To non-dimensionalise the system, consider the units of state variable $\hat{V}$ and the parameter $\alpha$: $[\hat{V}] =$\# virions, $[\alpha] = $\#virions per cell. As such, $\hat{V}$ can be re-scaled to give:

\begin{equation}
V=\frac{v}{\alpha}, \ \ \ \  [V] = \frac{[v]}{[\alpha]} = \frac{\#\mbox{virions}}{\#\mbox{virions per cell}} = \# \mbox{cells},
\end{equation}
\noindent which has no units. Substituting this re-scaled variable into the model gives:

\begin{equation}
\begin{split}\nonumber
\frac{dU}{d\tau} &= r\ln\left(\frac{K}{U}\right)U - \frac{ \beta\alpha UV}{U+I}\\
\frac{dI}{d\tau} &= \frac{\beta\alpha  UV}{U+I}-d_I I\\
\frac{dV}{d\tau} &= d_I I -d_V V
\end{split}
\end{equation}

\noindent To eliminate time, $\beta$ is used to rescale $\tau$. The units of $\beta$ are

\begin{equation}
[\beta] = \frac{1}{\#\mbox{virions}}\frac{1}{\mbox{per unit time}},
\end{equation}

\noindent so for $\hat{\beta}=\beta \alpha$, the units of $\hat{\beta}$ would be $[\hat{\beta}]=1/\text{time}$. Scaling time by $\hat{\beta}$ gives $[t] = [\hat{\beta}][\tau]$ which is dimensionless and the system of equations can now be rewritten for the independent dimensionless variable $t$:

\begin{align}\nonumber
\frac{dU}{dt} &= \frac{r}{\hat{\beta}}\ln\left(\frac{K}{U}\right)U - \frac{UV}{U+I},\\
\frac{dI}{dt} &= \frac {UV}{U+I}-\frac{d_I}{\hat{\beta}}I,\\
\frac{dV}{dt} &= \frac{d_I}{\hat{\beta}}I -\frac{d_V}{\hat{\beta}}V,
\end{align}

\noindent which is equivalent to Eqs.~(\ref{E4})-(\ref{E6}), where
\begin{equation}\nonumber
m = \frac{r}{\hat{\beta}}, \ \ \xi = \frac{d_I}{\hat{\beta}}, \ \ \ \gamma = \frac{d_V}{\hat{\beta}}. 
\end{equation}

\noindent In this article, we have chosen to consider the model in this form; however, it is also possible to eliminate the dimensions for $U, I$ and $V$. 

Currently, $[U]$ = $[I]$ = $[V]$ = \# cells, and $[K] =$ \#cells. Denote the dimensionless forms by 

\begin{equation}\nonumber
[\tilde{U}] =\frac{[U]}{[K]} = \frac{\#\mbox{cells}}{\#\mbox{cells}}, \ \ [\tilde{I}] =\frac{[I]}{[K]} = \frac{\#\mbox{cells}}{\#\mbox{cells}}, \ \ [\tilde{V}] =\frac{[V]}{[K]} = \frac{\#\mbox{cells}}{\#\mbox{cells}}
\end{equation}

\noindent So $\tilde{U}, \tilde{I}, \tilde{V}$ are dimensionless. By dividing all terms by $K$, the model becomes:

\begin{align}\nonumber
\frac{d\tilde{U}}{dt} &= r\ln\left(\frac{1}{\tilde{U}}\right)\tilde{U} - \frac{ \beta\alpha \tilde{U}V}{U+I}\\
\frac{d\tilde{I}}{dt} &= \frac{\beta\alpha \tilde{ U}V}{U+I}-d_I\tilde{ I}\\
\frac{d\tilde{V}}{dt} &= d_I\tilde{ I} -d_V \tilde{V}
\end{align}

\noindent and since

\begin{equation}
\frac{V}{U+I}=\frac{K\tilde{V}}{K\tilde{U}+K\tilde{I}}= \frac{\tilde{V}}{\tilde{U}+\tilde{I}}
\end{equation}

\noindent the full non-dimensional system of equations is finally given by:

\begin{equation}
\begin{split}\nonumber
\frac{d\tilde{U}}{dt} &= m\ln\left(\frac{1}{\tilde{U}}\right)\tilde{U} - \frac{\tilde{U}\tilde{V}}{\tilde{U}+\tilde{I}}\\
\frac{d\tilde{I}}{dt} &= \frac{\tilde{ U}\tilde{V}}{\tilde{U}+\tilde{I}}-\xi\tilde{ I}\\
\frac{d\tilde{V}}{dt} &= \xi\tilde{ I} -\gamma \tilde{V}
\end{split}
\end{equation}

\newpage

\bibliographystyle{plain}
\bibliography{references}

\end{document}